\documentclass[aos,MSNbibl,citesort,seceqn,dvips]{arximspdf}
\usepackage{multirow}
\usepackage{graphicx}

\doi{10.1214/11-AOS882}
\volume{39}
\issue{3}
\pubyear{2011}
\firstpage{1720}
\lastpage{1747}

\makeatletter

\newcommand{\Var}{\operatorname{var}}
\def\a{\alpha}
\def\approx{^{\mathrm{approx}}}
\def\be{\beta}
\def\biggmi{\mid}
\def\Bigmi{\mid}
\def\bep{{\bar\ep}}
\def\bga{{\bar\ga}}
\def\bX{{\bar X}}
\def\bY{{\bar Y}}
\def\cB{{\cal B}}
\def\cI{{\cal I}}
\def\cJ{{\cal J}}
\def\cX{{\cal X}}
\def\de{\delta}
\def\ep{\epsilon}
\def\error{\mathrm{error}}
\def\ga{\gamma}
\def\hbe{{\hat\be}}
\def\hg{{\hat g}}
\def\hga{\hat{b}}
\def\half{^{1/2}}
\def\inti{\int_{\cI}}
\def\ij{_{ij}}
\def\J{^{\cJ}}
\def\mi{\mid}
\def\mo{^{-1}}
\def\one{^{[1]}}
\def\otd{(1/3)}
\def\ra{\to}
\def\rai{\to\infty}
\def\si{\sigma}
\def\sumcj{\sum_{j\in\cJ(\be)}}
\def\sumi{\sum_i }
\def\sumij{\mathop{\sum\sum}_{i,j\dvtx i\neq j}}
\def\suminj{\sum_{i\dvtx i\neq j}}
\def\sumion{\sum_{i=1}^n }
\def\sumj{\sum_j }
\def\sumjon{\sum_{j=1}^n }
\def\ta{{\zeta}}
\def\tb{{\vartheta}}
\def\th{\theta}
\def\z{^0}
\newcommand{\cal}{\mathcal}
\newcommand{\eqref}[1]{(\ref{#1})}
\renewcommand{\epsilon}{\varepsilon}
\newtheorem{theorem}{Theorem}[section]
\makeatother

\begin{document}
\begin{frontmatter}

\title{Single and multiple index functional regression models with nonparametric link}
\runtitle{Functional regression with nonparametric link}

\begin{aug}
\author[A]{\fnms{Dong} \snm{Chen}\corref{}\ead[label=e1]{dchen@wald.ucdavis.edu}},
\author[B]{\fnms{Peter} \snm{Hall}\thanksref{t1}\ead[label=e2]{halpstat@ms.unimelb.edu.au}}
\and
\author[A]{\fnms{Hans-Georg} \snm{M\"{u}ller}\thanksref{t2}\ead[label=e3]{mueller@wald.ucdavis.edu}}
\runauthor{D. Chen, P. Hall and H.-G. M\"{u}ller}
\affiliation{University of California, Davis, University of California,
Davis and University~of~Melbourne, and University of California, Davis}
\address[A]{D. Chen\\
H.-G. M\"{u}ller\\
Department of Statistics\\
University of California, Davis\\
Davis, California 95616\\
USA\\
\printead{e1}\\
\phantom{\textsc{E-mail:}\ }\printead*{e3}} 
\address[B]{P. Hall\\
Department of Mathematics\\ \quad and Statistics\\
University of Melbourne\\
Parkville, VIC 3010\\
Australia\\
\printead{e2}}
\end{aug}
\thankstext{t1}{Supported in part by an Australian Research Council
Fellowship.}
\thankstext{t2}{Supported in part by National Science Foundation Grant
DMS-08-06199.}

\received{\smonth{7} \syear{2010}}
\revised{\smonth{2} \syear{2011}}

%
\begin{abstract}
Fully nonparametric methods for regression from functional data have
poor accuracy from a statistical viewpoint, reflecting the fact that
their convergence rates are slower than nonparametric rates for the
estimation of high-dimensional functions. This difficulty has led to
an emphasis on the so-called functional linear model, which is much
more flexible than common linear models in finite dimension, but
nevertheless imposes structural constraints on the relationship
between predictors and responses. Recent advances have extended the
linear approach by using it in conjunction with link functions, and
by considering multiple indices, but the flexibility of this
technique is still limited. For example, the link may be modeled
parametrically or on a grid only, or may be constrained by an
assumption such as monotonicity; multiple indices have been modeled
by making finite-dimensional assumptions. In this paper we introduce
a new technique for estimating the link function nonparametrically,
and we suggest an approach to multi-index modeling using adaptively
defined linear projections of functional data. We show that our
methods enable prediction with polynomial convergence rates. The
finite sample performance of our methods is studied in simulations,
and is illustrated by an application to a functional regression problem.
\end{abstract}

%
\begin{keyword}[class=AMS]
\kwd[Primary ]{62G05}
\kwd{62G08}.
\end{keyword}
\begin{keyword}
\kwd{Functional data analysis}
\kwd{generalized functional linear model}
\kwd{prediction}
\kwd{smoothing}.
\end{keyword}

\end{frontmatter}

\section{Introduction}\label{sec1}

When explanatory variables are functions, rather than
vectors, the problems of nonparametric regression and prediction
are intrinsically difficult from a statistical viewpoint. In
particular, convergence rates can be slower than the inverse of
any polynomial in sample size, and so relatively large samples may
be needed in order to ensure adequate performance. Fully
nonparametric methods have been studied recently in
functional data regression and related problems (see, e.g.,
\cite{9,7} and~\cite{14}). The slow convergence rates
associated with these unstructured nonparametric approaches
provide motivation for seeking nonparametric approaches that
exploit a greater amount of structure in the data and are
therefore expected to have better properties from a statistical
perspective.

Advances in this direction include those made in~\cite{18,19,17,11} and~\cite{13}, where both parametric and nonparametric
link functions were introduced in order to connect the response to a linear
functional model in the explanatory variables. However, the flexibility of
available link-function models is still rather limited. For example,
although nonparametric link functions were considered in~\cite{13}, the
approaches considered there are restricted by the assumption of monotonicity,
where the corresponding ``Generalized Functional Linear Model''
approach is
based on a semiparametric quasi-likelihood based estimating equation, which
includes known or unknown link and variance functions. In contrast, we are
aiming here at models with one or several nonparametric link functions,
ignoring possible heteroscedasticity of the errors. Our approach
provides an
alternative to the related methods in~\cite{ait08}, where single-index
functional regression models with general nonparametric link functions are
considered that may be chosen nonmonotonically and without shape
constraints. The main differences are that our methodology includes the
multi-index case, does not anchor the true parameter on a prespecified
sieve, and that we provide a detailed theoretical analysis of a direct
kernel-based estimation scheme that culminates in a convergence result that
establishes a polynomial rate of convergence.

Beyond demonstrating that our approach enables prediction with
polynomial accuracy, we also include generalizations to
iteratively fitted multiple index models, founded on a sequence of
linear regressions. Here we borrow ideas from dimension reduction
in models that involve high-dimensional, but not functional, data.
When the link function is nonparametric, the intercept term in
functional linear regression loses its relevance because it is
incorporated into the link. The slope function is still
potentially of interest, but the viewpoint taken in this paper is
predominately one of prediction rather than slope estimation, and
in particular our theory focuses directly on the prediction
problem. We refer to the papers by~\cite{21,card} and
\cite{cramb} for further discussion of these objectives in the
context of the functional linear model.

We introduce our model and estimation methodology in Section~\ref{sec2}.
Theoretical results regarding the polynomial convergence rate are
discussed in Section~\ref{sec3}, while algorithmic details are described in
Section~\ref{sec4}, which also includes an illustration of the proposed
methods with an application to spectral data. Simulation results
are reported in Section~\ref{sec5}. Detailed assumptions and proofs can be
found in the \hyperref[appm]{Appendix}.

\section{Model and methodology}\label{sec2}

\subsection{Model}\label{sec2.1}

Suppose we observe data pairs $(X_1,Y_1),\ldots,(X_n,Y_n)$,
independent and
identically distributed as $(X,Y)$, where $X$ is a random function
defined on
a compact interval $\cal{I}$ and $Y$ is a scalar. We anticipate that $(X,Y)$
is generated as
%
\begin{equation}\label{2.1} Y=g(X)+\error ,
\end{equation}
where $g$ is a smooth
functional and the error has zero mean, finite variance and is uncorrelated
with $X$. The model at \eqref{2.1} admits many interpretations and
generalizations, where, for instance, $X$ is a multivariate rather than
univariate function. For example, $X$ might be $(Z,Z')$, where $Z$ is a
univariate function and $Z'$ its derivative. To simplify the developments,
we shall focus on problems where $X$ is a univariate function of a single
variable. Models and methodology in more general settings are readily
developed from the single-variable case. Our approach is described in detail
for situations where the trajectories of functional predictors can be assumed
to be fully observed, for example, due to smoothness such as for the Tecator
data which we analyze with the proposed methods in Section~\ref{sec4.2}; it can be
extended to cases with densely and regularly measured trajectories, where
measurements may be subject to i.i.d. noise with finite fourth-order moments.
This extension requires sufficiently dense measurement designs, such that
smoothness assumptions coupled with suitable smoothing methods lead to
sufficiently fast uniform rates of convergence when pre-smoothing the
data to
generate smooth trajectories. Such an extension will not be feasible for
functional data for which only sparse and noisy measurements are available.

The case where $g$, in \eqref{2.1}, is a general functional, even a very
smooth functional, can have serious drawbacks from the viewpoint of practical
function estimation, since the problem of estimating such a $g$ is inherently
difficult from a statistical viewpoint. In particular, convergence
rates of
estimators in this case are generally slower than the inverse of any
polynomial in sample size. Therefore, unless the data set is very
large, it
can be particularly difficult to estimate $g$ effectively. In this respect
the commonly assumed functional linear model, where $g(x)=\a+\inti\be
 x$,
$\a$ is a scalar and $\be$ is a regression parameter function, offers
substantial advantages, for example, polynomial convergence rates and
even, on
occasion, root-$n$ consistency. However, the linear-model assumption is often
too restrictive in practical applications.

An alternative approach is to place the linear model inside a link
function, for example, defining
%
\begin{equation}\label{2.2}
g(x)=g_1 \biggl(\a+\inti\be x \biggr) ,
\end{equation}
although this, too, is
restrictive unless we select the link in a very adaptive manner.
We suggest choosing\vadjust{\goodbreak} the link function $g_1$ nonparametrically. In
this case the intercept parameter, $\a$, in \eqref{2.2} is
superfluous; it can be replaced by zero, and its effect
incorporated into $g_1$. Therefore we actually fit the model
%
\begin{equation}\label{2.3} g(x)=g_1 \biggl(\inti\be x \biggr) ,
\end{equation}
where $g_1$
is subject only to smoothness conditions, and to ensure
identifiability, we require a condition on the ``scale'' of
$\beta$, which we choose as $\inti\be^2=1$. The sign of $\beta$
can be determined arbitrarily.

\subsection{Methodology}\label{sec2.2}

We estimate the parameter function $\be$ and the link function
$g_1$ in the model at \eqref{2.3}, using least squares in
conjunction with local-constant or local-linear smoothing as
follows. To obtain $g_1$, we will use a scatterplot smoother
which we implement as local-constant or local-linear weighted
least squares smoothing. Given a parameter function $\be$, the
scatterplot targeting the nonparametric regression
$g_1(z)=E(Y|\inti\be  X=z)$ consists of the data pairs
$(\inti\be  X_i, Y_i)_{i=1,\ldots,n}.$ Omitting the predictor
$X_j$ when predicting the response at $\inti\be  X_j$, averaging
least squares smoothers constructed for predicting at the observed
predictor levels $X_j$ are then obtained by choosing intercept
parameters $\ta_j$ and slope parameters $\tb_j$ to minimize
%
\begin{eqnarray}\label{2.4}
&\displaystyle\sumij
(Y_i-\ta_j)^2 K \biggl\{h\mo\inti\be (X_i-X_j) \biggr\} \quad\mbox{or}&
\nonumber
\\[-8pt]
\\[-8pt]
&\displaystyle\sumij \biggl\{Y_i- \biggl(\ta_j+\tb_j\inti\be X_i\biggr ) \biggr\}
^2
K \biggl\{h\mo\inti\be (X_i-X_j) \biggr\} ,&
\nonumber
\end{eqnarray}
in the local-constant and local-linear cases, respectively, where
$K$ is a kernel function and $h$ is a bandwidth.

Defining $K\ij=K\{h^{-1}\inti\be (X_i-X_j)\}$, $\bX_j=(\suminj
X_i K\ij)/\suminj K\ij$ and $\bY_j=(\suminj Y_i K\ij)/\suminj
K\ij$, one finds that the minimia of \eqref{2.4}, for any given $\be$,
are
%
\begin{equation}\label{2.5}
  \sumij(Y_i-\bY_j)^2 K\ij\quad\mbox{or}\quad
\sumij \biggl\{Y_i-\bY_j
-\hat{\tb}_j \inti\be (X_i-\bX_j) \biggr\}^2 K\ij.\hspace*{-35pt}
\end{equation}
The
minimizers $\hat{\ta}_j$ are given by
$\hat{\ta}_j=\hat{\ta}_j(\be)=\bY_j$ in the local-constant case
and in the local-linear case minimizers $\hat{\ta}_j$ and
$\hat{\tb}_j$ are given by
%
\begin{eqnarray}\label{2.6}
 \hat{\ta}_j(\be)&=&\bY_j-\hat{\tb}_j(\be)\inti\be\bX_j,\nonumber
\\[-8pt]
\\[-8pt]
 \hat{\tb}_j(\be)&=&\frac{\suminj\{\inti\be(X_i-\bX_j)\}
(Y_i-\bY_j)
K\ij}{\suminj\{\inti\be (X_i-\bX_j)\}^2 K\ij} , \qquad 1\leq
j\leq n.
\nonumber
\end{eqnarray}
Summarizing, the criteria at \eqref{2.4} are based on
averaging local-constant and local-linear fits to
$g_1(\inti\be x)$ at $x=X_j$, averaging over\vadjust{\goodbreak} $X_j$, where the
respective fits are computed from the data $X_1,\ldots,X_n$,
excluding $X_j$. The resulting approximations to $g_1(\inti\be
 X_j)$, for a given $\be$, are $\bY_j$ and
$\bY_j+\hat{\tb}_j(\be) \inti\be (X_j-\bX_j)$, respectively,
with $\hat{\tb}_j(\be)$ given by \eqref{2.6}.

It remains to specify our final estimates. We estimate $\be$ by conventional
least squares, aiming to minimize the sum of squared differences between
$Y_j$ and the just-mentioned approximations:
%
\begin{eqnarray}\label{2.7}
 S(\be)&=&\sumjon\bigl(Y_j-\bY_j(\be)\bigr)^2\quad\mbox{or} \nonumber
\\[-8pt]
\\[-8pt]
 S(\be)&=&\sumjon \biggl\{Y_j-\bY_j(\be)
-\hat{\tb}_j(\be) \inti\be (X_j-\bX_j) \biggr\}^{ 2} ,
\nonumber
\end{eqnarray}
subject to
$\inti\be^2=1$ and with $\hat{\tb}_j(\be)$ as in \eqref{2.6}. This
problem is most
conveniently solved by expanding $\be=\sum_{1\leq k\leq r} b_k \psi_k$,
where $\psi_1,\psi_2,\ldots$ is an orthonormal basis and $r$ denotes a
frequency cut-off, choosing the generalized Fourier coefficients $b_k$ to
minimize $S(\be)$. This gives estimators $\hga_1,\ldots,\hga_r$ of
$b_1,\ldots,b_r$, respectively, and from those we may compute our estimator
of $\be$:
%
\begin{equation}\label{2.8} \hbe=\sum_{k=1}^r \hga_k \psi_k
\qquad
\mbox{subject to } \sum_{k=1}^r \hga_k^2=1.
\end{equation}

The basis $\psi_1,\psi_2,\ldots$ can be chosen as a fixed basis such
as one
of various orthonormal polynomial systems or the Fourier basis, or
could be
another sequence altogether, chosen for computational convenience. We note
that it does not matter for the validity of our results whether the basis
functions are fixed or random. Therefore the basis can be chosen as estimated
eigenfunction basis, as long as the estimated eigenfunctions are orthonormal.
We note that irrespective of how it is constructed, the selected basis needs
to be such that condition \eqref{3.4} below is satisfied for the generalized
Fourier coefficients of $\beta$, while the additionally needed conditions
\eqref{3.5}, \eqref{3.6} depend only on properties of $\beta$ and $X$
but not
on the choice of the basis. The condition at \eqref{3.4} requires a
polynomial decay rate (of arbitrary order) for the tail sums of the Fourier
coefficients of $\beta$, which is slightly stronger than the
convergence of
the tail sums to 0 that is implied by the square integrability of
$\beta$.
Since we do not assume prior knowledge about $\beta$, no particular
basis is
preferable in this regard a priori. In any case, the theory applies if
\eqref{3.4} holds for the selected basis. In practice, one would choose
a basis based on
how well the representation of $\beta$ works in typical applications. We
found the choice of estimated eigenfunctions for representing $\beta$
particularly convenient for our applications and our implementation is
therefore using this basis.

We note that the criteria at \eqref{2.4} are not directly comparable with
those at \eqref{2.7}, not least because in \eqref{2.4} we are fitting $g_1$
locally and in \eqref{2.7} we are fitting $\be$ globally. Reflecting these
two different contexts, each residual squared error in both criteria in
\eqref{2.4} has a local kernel weight, whereas each residual squared
error in the criteria in \eqref{2.7} has a constant weight.

Having computed $\hbe$, we estimate the univariate function $g_1(u)$ by
conventional local-constant or local-linear regression on the pairs
$(\inti\hbe X_i,Y_i)$, for $1\leq i\leq n$. In particular, in the
local-constant case we take
%
\begin{equation}\label{2.9} \hg_1(u)= \Biggl\{\sumion
Y_i K_i(u) \Biggr\}\bigg /  \Biggl\{\sumion K_i(u) \Biggr\} ,
\end{equation}
where
$K_i(u)=K\{h\mo (\inti\hbe X_i-u)\}$; in the local-linear setting
we choose
$\ta=\hat{\ta}$ and $\tb=\hat{\tb}$, both of which can also be
viewed as
functions of $u$, to minimize
$\sumi\{Y_i-(\ta+\tb\inti\hbe X_i)\}^2 K_i(u)$, and then put
$\hg_1(u)=\hat{\ta}+\hat{\tb} u$. Several aspects of this algorithm
can be
modified to improve its performance. For example, noting that the ratio on
the right-hand side of \eqref{2.6} will likely be unstable if the denominator
is based on a relatively small number of terms, we might restrict the sum
over $j$ in either formula in \eqref{2.5} to values of that index for which
$\suminj K\ij\geq{\lambda }$, where ${\lambda }>0$ denotes a
sufficiently large threshold,
and repeat this restriction in the case of \eqref{2.7}. Problems caused
by a
too-small denominator can be especially serious in the case of functional
data, since sample sizes there are typically relatively small.

If we take the view that the problem of interest is that of estimating $g$
for the purpose of prediction, and that estimating $\be$ and $g_1$ in their
own right is of relatively minor interest, then standard
cross-validation can
be used to choose simultaneously the smoothing parameters $h$ and $r$. In
Section~\ref{sec3} we adopt the perspective of prediction, and show that in that
context the estimator $\hg$ approximates $g$ at a rate that is polynomially
fast as a function of sample size.

\subsection{Multiple index models}\label{sec2.3}

The model at \eqref{2.3} can be generalized by taking $g_1$ to be
a $p$-variate function
%
\begin{equation}\label{2.10}
g(x)=g_1 \biggl(\inti\be_1 x,\ldots,\inti\be_p x \biggr), \qquad
\inti\be_j^2=1 \mbox{ for } 1\leq j\leq p.
\end{equation}
However,
given the relatively small sample sizes often encountered in
functional data analysis, focusing on the function at (\ref{2.10}), with
$p\geq2$, will often lead to estimators with high variability. An
alternative, $p$-component functional multiple index model,
such as
%
\begin{equation}\label{2.11} g(x)=g_1 \biggl(\inti\be_1 x
\biggr)+\cdots
+g_p \biggl(\inti\be_p x \biggr) ,
\end{equation}
is arguably more attractive.
This class of models has been considered by~\cite{23}, who
referred to them as ``Functional Adaptive Models.'' The approach of
James and Silverman was restricted to the parametric case by
requiring the functional predictors $x_i$ as well as the link
functions $g_j$ to be elements of a finite-dimensional spline
space, excluding nonparametric (infinite-dimensional) link and
predictor functions. Such a fully parametric framework allows the
use of a likelihood-based approach to fitting these models,
establishing identifiability by extending previous results for the
vector case~\cite{20}.

Since our main goal is prediction and not model
identification, we are not primarily concerned with
identifiability issues and do not emphasize specific
identifiability conditions for the models we consider. The models
at \eqref{2.10} and \eqref{2.11} in fact are not identifiable
without further restrictions. To appreciate why, note that the
order of the components on the right-hand side of \eqref{2.10}, or
of the functions on the right-hand side of \eqref{2.11}, could be
permuted without affecting the model. This problem does not arise
for conventional multivariate or additive models, where the
arguments of the functions are predetermined as the components of
the explanatory variable~$x$. While this difficulty can be
overcome in a variety of ways, using a recursive additive model is
attractive on both statistical and computational grounds. We now
give background for that approach.

It is not uncommon in statistics to pragmatically alter a difficult problem
to one that is simpler. Indeed, the introduction of additive models is
typically motivated in that manner. Thus, we could generalize the
problem of
estimating a link function~$g$, and a slope function $\be$, in \eqref{2.1},
subject only to smoothness conditions, to that of estimating the
intrinsically simpler functions defined at \eqref{2.11}.
Alternatively, and
more appropriately from the perspective of general inference, we would seek
to estimate $g$ in \eqref{2.10} not because we felt that those
functions were
identical to $g$ in \eqref{2.1}, but because they were relatively accessible
approximations to $g$. Taking this view of the problem of estimating, or
rather, approximating, the function $g$ in \eqref{2.1}, and accepting that
the $p$-additive function at \eqref{2.11} is more likely to be
practicable in
functional data analysis than the $p$-variate function at (\ref{2.10}), we suggest
fitting the $g$ in \eqref{2.11} recursively, for steadily increasing values
of $p$. This ``backfitting'' approach borrows an idea from projection
pursuit regression, to use recursive, low-dimensional, projection-based
approximations.

In particular, taking $g_1^0=g\z$ where $g\z$ denotes the true value
of $g$ at \eqref{2.1}, we choose the function $g_1$ of a single
variable, and the function $\be_1$, to minimize, in the case $j=1$,
the expected value
%
\begin{equation}\label{2.12}
E \biggl\{g_j^0(X)-g_j \biggl(\inti\be_j X \biggr) \biggr\}^2 \qquad
\mbox{subject to } \inti\be_j^2=1.
\end{equation}
More generally, if we
have calculated $\be_{j-1}$ and $g_{j-1}$, and previously defined
$g_{j-1}^0(x)$, then we may define $g_j^0$ by
$g_j^0(x)=g_{j-1}^0(x)-g_{j-1}(\inti\be_{j-1} X)$ and choose $g_j$
and $\be_j$ to minimize the quantity at \eqref{2.12}.

In the next section we shall show how to calculate estimators
$\hg_j$ and $\hbe_j$ of $g_j$ and $\be_j$, respectively, for
$j\geq1$. Note that we do not claim to consistently estimate $g$,
in \eqref{2.1}, unless that function has exactly the form at
\eqref{2.3} (in which case our estimator is $\hg=\hg_1$, defined in
Section~\ref{sec2.2}). Instead we suggest developing consistent estimators of
successive approximations to $g(x)$, that is, of
%
\begin{eqnarray}\label{2.13}
&&  g_1 \biggl(\inti\be_1 x \biggr) ,
 g_1 \biggl(\inti\be_1 x \biggr)+g_2 \biggl(\inti\be_2 x \biggr) ,\nonumber
\\[-8pt]
\\[-8pt]
 &&g_1 \biggl(\inti\be_1 x \biggr)+g_2 \biggl(\inti\be_2 x \biggr)
+g_3 \biggl(\inti\be_3 x \biggr) ,\ldots .
\nonumber
\end{eqnarray}

\subsection{Estimation in functional multiple index models}\label{sec2.4}

Here we generalize the methodology in Section~\ref{sec2.2} so that it permits
estimation of the functions $g_1,g_2,\ldots$ in \eqref{2.12}. Assume
we are fitting a $p$-index model. The recursive fitting procedure
means once we have estimators $\hbe_j$ and $\hg_j$, for $1\leq j\leq
k-1<p$, of the functions $\be_j$ and $g_j$ defined in the paragraph
containing \eqref{2.12}, we take
$Y_i(k)=Y_i-\hg_1(X_i)-\cdots-\hg_{k-1}(X_i)$, and use the
methodology in Section~\ref{sec2.2} but with $Y_i(k)$ replacing $Y_i$,
obtaining an estimator $\hbe$, on this occasion actually an
estimator $\hbe_k$ of $\be_k$, and an estimator $\hg$, which is
really an estimator $\hg_k$ of $g_k$. The quantity
$\hg_1(\inti\hbe_1 x)+\cdots
+\hg_p(\inti\hbe_p x)$ is our estimate of the $p$-index
model from the recursive fitting procedure.

A further refinement that leads to smaller prediction errors is
backfitting, which uses the recursive fits described above as a
starting point. Once these fits are obtained, further updates are
obtained iteratively by revisiting and updating one index after
another, presuming that the remaining $p-1$ indexes are fixed. The
iterative updating of individual indices is itself iterated until
indices change only little. This is implemented in a similar way as
described in~\cite{20} for a traditional multiple index model with
monotone link functions. Denoting the estimates obtained from the
initial recursive fitting procedure by
$\hg^0_1(\inti\hbe^0_1 x)+\cdots
+\hg^0_p(\inti\hbe^0_p x)$, then for the $d$th
iteration, iterating also through the increasing sequence
$k=1,2,\ldots,p$, one uses
%
\begin{equation}\label{2.14} Y^d_i(k)=Y_i-\sum_{j<k}
\hg^{d}_j \biggl(\inti\hbe^{d}_j X_i \biggr) - \sum_{j>k}
\hg^{d-1}_j \biggl(\inti\hbe^{d-1}_j X_i \biggr)
\end{equation}
to replace $Y_i$
for fitting $\hg^d_k(\inti\hbe^d_k x)$. The iterative
backfitting procedure is stopped once the relative differences
between $\hbe_1^{d-1}$ and $\hbe_1^{d}$ fall below a prespecified
threshold or a maximum number of iterations is reached.

\section{Polynomial convergence rate}\label{sec3}

The main result in this section establishes that, if the linear model is
linked to the response variable as in \eqref{2.3}, if a H\"{o}lder smoothness\vadjust{\goodbreak}
condition on the link function $g_1$ is assumed, and if we ask of the
generalized Fourier expansion $\be=\sum_{k\geq1} b_k \psi_k$ that it
converges polynomially fast at a sufficiently rapid rate, then the predictor
$\hg$ converges to $g$ at a polynomial rate. That property
distinguishes the
approach suggested in this paper from fully nonparametric methods that impose
only smoothness conditions on the function $g$, in \eqref{2.1}, but
have much
slower convergence rates for the predictor. We give explicit theory
only in
the local-constant case, since, as argued at the end of Section~\ref{sec2.2}, that
approach is particularly appropriate when dealing with functional data. The
local-linear setting can be treated similarly.

We assume that independent and identically distributed data pairs
$(X_i,Y_i)$ are generated by the model discussed in Section~\ref{sec2}:
%
\begin{equation}\label{3.1}
\begin{tabular}{@{}p{297pt}@{}}
$Y_i=g(X_i)+\ep_i$, where the $X_i$'s are
square-integrable
random functions supported on the
compact interval $\cI$,
$g$ is a real-valued functional given by $g(x)=g_1(\inti\be\z
x)$,
$g_1$ is a
real-valued\vspace*{-1pt} function of a single variable, $\be\z$ enjoys
the property
$\inti{\be\z}^2=1$ and
denotes the true value of
the square-integrable function $\be$, and
the errors $\ep_i$ are
independent of the $X_i$'s and of one another, and have
zero mean.
\end{tabular}
\end{equation}
The only assumption we make of $g_1$ is that it is bounded and
smooth:
%
\begin{equation}\label{3.2}
\begin{tabular}{@{}p{297pt}@{}}
$g_1$ is bounded and satisfies a H\"{o}lder continuity
condition: $|g_1(u)-g_1(v)|\leq D_1 |u-v|^{a_1}$ for all
$u$ and $v$, where $a_1,D_1>0$.
\end{tabular}
\end{equation}
The assumption that $g_1$ is bounded can be relaxed. For example,
if the functions $X_i$ are bounded with probability~1, then
$\inti\be\z X_i$ is uniformly bounded, and so the distribution of
the response variables $Y_i$ depends only on the values that $g_1$
takes on a particular compact interval. We can extend $g_1$ from
that interval to the whole real line in such a way that the extended
version of $g_1$ is bounded and has a bounded derivative. More
generally, if $\sup_{1\leq i\leq n} \|X_i\|$ grows at rate
$O(n^\eta)$, for all $\eta>0$, where $\|X\|$ denotes the $L_2$ norm
of $X$ (e.g., this condition holds if $X$ is a Gaussian
process), and if $\sup_{|x|\leq u} |g_1(x)|$ grows at no faster
than a polynomial rate as $u$ diverges, then only minor
modifications of our proof of the theorem are required to establish
Theorem~\ref{theo3.1}.

Let $X$ have the common distribution of the random functions $X_i$
in the model at \eqref{3.1}. We ask that $\|X\|$ have at least a
small, fractional moment, and that all moments of the error
distribution be finite. In particular:
%
\begin{equation} \label{3.3}
\begin{tabular}{@{}p{297pt}@{}}
$E(\|X\|^\eta)<\infty$ for some $\eta>0$, and
$E(|\ep|^m)\leq(D_2 m)^{a_2m}$
for all integers $m\geq1$, where $a_2,D_2$ denote
positive constants.
\end{tabular}
\end{equation}
The condition $E|\ep|^m\leq(D_2 m)^{a_2m}$ is satisfied by
distributions the tails of which decrease at rate at least
$\exp(-C_1 x^{C_2})$, for constants $C_1,C_2>0$, provided we choose\vadjust{\goodbreak}
$a_2>1/C_2$. In particular, the condition is satisfied by
exponential and Gaussian distributions, and also, in the case
$C_2<1$, by many distributions that do not have finite moment
generating functions.

Write $f( \cdot\mi\be)$ for the probability density of $\inti\be
X$. Given an orthonormal basis $\psi_1,\psi_2,\ldots$ for the
class $L_2(\cI)$ of square-integrable functions on $\cI$, express
a general function $\be\in L_2(\cI)$ with $\inti\be^2=1$ as
$\be=\sum_{k\geq1} b_k \psi_k$, where $\sum_{k\geq1} b_k^2=1$.
For constants $a_3,a_4,B,D_3,D_4,D_5>0$, we shall assume that:
%
\begin{eqnarray}\label{3.4}
&\displaystyle\sum_{k=r+1}^\infty b_k^2\leq D_3 (1+r)^{-B} \qquad\mbox{for all }
r\geq1 ,&\\\label{3.5}
&\displaystyle\sup_{\be\in\cB ; x} f(x\mi\be)<\infty ,&\\\label{3.6}
&\displaystyle\sup_{\be\in\cB} P \biggl\{f \biggl(\inti\be
X-u  | \be \biggr)\leq D_4 \de^{a_3} \mbox{ for all }
|u|\leq\de \biggr\}\leq D_5 \de^{a_4} ,&
\end{eqnarray}
where \eqref{3.6} holds for all sufficiently small $\de>0$. Condition
\eqref{3.4} is standard; it asks that the generalized Fourier
coefficients of $\be$ decay at least polynomially fast, in a weak
sense. To appreciate the motivation for \eqref{3.5} and \eqref{3.6}, observe
that if $X$ is a Gaussian process for which the covariance
operator has eigenvalues $\th_1\geq\th_2\geq\cdots\geq0$ and
respective eigenfunctions $\phi_1,\phi_2,\ldots,$ then
$f( \cdot\mi\be)$ is the N$(a,\varsigma^2)$ density, where
$a=\inti\be E(X)$, $\varsigma^2=\sum_{k\geq1} \th_k b_k^2$ and
$b_k=\inti\be \phi_k$. Then \eqref{3.6} is obtained by using
well-known tail bounds for the Gaussian distribution function
$\Phi$ with standard Gaussian density $\phi$. It follows that
\eqref{3.5} and \eqref{3.6} hold whenever $0<a_4\leq a_3<\infty$ and
$\cB$ is
a class of functions $\be$ for which $\sum_{k\geq1} \th_k b_k^2$
is bounded away from zero and infinity, and for which
$\sum_{k\geq1}b_k^2=1$. Our use of the principal component basis
in this example serves only to show the reasonableness of
conditions \eqref{3.5} and \eqref{3.6}, which of course do not depend on
choice of basis. It does not imply that the basis
$\psi_1,\psi_2,\ldots$ should be identical to
$\phi_1,\phi_2,\ldots.$

Of the kernel $K$ and bandwidth $h$ we ask that:
%
\begin{equation} \label{3.7}
\begin{tabular}{@{}p{297pt}@{}}
$K$ is nonnegative and symmetric, has support equal to
a compact interval, decreases to zero as a polynomial
at the ends of its support, and has a bounded derivative;
and $h\sim
D_6 n^{-C}$ as $n\rai$, where $C,D_6>0$.
\end{tabular}
\end{equation}

Define $\hbe$ to be the minimizer of $S(\be)=\sumj(Y_j-\bY_j)^2$
[the first quantity in \eqref{2.7}, corresponding to
local-constant estimation] over functions $\be=\sum_{1\leq k\leq
r} b_k \psi_k$, constrained by $\sum_{1\leq k\leq r} b_k^2=1$,
for which \eqref{3.5} and \eqref{3.6} hold and $\sup_{k\geq
1} |b_k|\leq D_7$, with $D_7>D_3$ [the latter as in \eqref{3.4}],
and where
%
\begin{equation}\label{3.8}
\begin{tabular}{@{}p{297pt}@{}}\quad
$r$ denotes the integer part of $D_8 n^D$, for constants
$D,D_8>0$.
\end{tabular}
\end{equation}
This is the procedure for constructing $\hg$ suggested in the
argument leading to \eqref{2.8}, in the local-constant case.

\begin{theorem}\label{theo3.1}
If \eqref{3.1}--\eqref{3.8} hold, if $B$ in \eqref{3.4} is
sufficiently large, and if $C$ and $D$ in \eqref{3.7} and
\eqref{3.8} are sufficiently small (all three constants depending
only on $a_1,\ldots,a_4$), then there exists a constant $c>0$ such
that, as $n\rai$,
%
\begin{equation}\label{3.9} n\mo \sumjon\{g(X_j)-\hg(X_j)\}^2
=O_p (n^{-c} ) .
\end{equation}
\end{theorem}

The proof is in the \hyperref[appm]{Appendix}. It is possible to extend
Theorem~\ref{theo3.1}
to the recursive additive model case formulated in Section~\ref{sec2.4},
although the argument there is significantly longer. As
explained earlier, for the case of Gaussian predictors $X$, the
choices $a_1=1, a_2=1, a_3=1, a_4=1$ are possible and then by
choosing the other constants judiciously, observing the various
constraints, one finds that one may obtain the rate of convergence
in (\ref{3.9}) for $c$ with $c<1/4$. We do not pursue here the
question of the optimality of this rate of convergence. An
assumption that has been made throughout is that the predictor
trajectories are fully observed. This is an idealized situation.
It is possible to weaken this assumption, assuming that the
trajectories are sampled on a dense grid of points so that
integrals such as those appearing in \eqref{2.12} can be closely
approximated.

\section{Algorithmic implementation and data illustration}\label{sec4}

\subsection{Description of the algorithm}\label{sec4.1}

\subsubsection*{Step 1. Estimating $\beta$}

We assumed that $h$, $r$ and the basis $\{\psi_1,\ldots,\psi_r\}$
(we used eigenbasis in our implementation) in \eqref{2.8} and
\eqref{2.9} were given. We set
$\hat{\beta}=\sum_{k=1}^r\hat{b}_k\psi_k$, and the coefficients
$\hat{b}_1,\ldots,\hat{b}_r$ were estimated by minimizing
\eqref{2.7}. Those $Y_j$ with $\sum_{i\dvtx i\neq j}K_{ij}<\lambda$
were dropped from the minimization (we chose $\lambda=0.1$).
Letting $\xi_{ik}=\int\psi_k x_i$ and writing $S(\beta)$ in
\eqref{2.7} in terms of $b_1,\ldots,b_r$,
%
\begin{equation}
S(b_1,\ldots,b_r)=\frac{1}{n}\sum_{j=1}^n \biggl(Y_j-\sum_{i\neq
j}w_{ij} Y_i \biggr)^2 \label{4.1}
\end{equation}
for the local-constant case,
where
\[w_{ij}(b_1,\ldots,b_r,h)=\frac{K (h^{-1}\sum_{k=1}^r
b_k (\xi_{ik}-\xi_{jk}) )}{\sum_{l\neq
j}K (h^{-1}\sum_{k=1}^r b_k (\xi_{lk}-\xi_{jk}) )}
\]
are
the terms related to $b_1,\ldots,b_r$. For the local-linear case,
$S(b_1,\ldots,b_r)$ is more complicated, with similar subsequent
steps.

We note that $(b_1,\ldots,b_r,h)$ are not identifiable without
constraints, since
$w_{ij}(b_1,\ldots,b_r,h)=w_{ij}(cb_1,\ldots,cb_r,ch)$ for any
constant $c$. Meanwhile, if $K$ is symmetric,\vadjust{\goodbreak}
$w_{ij}(b_1,\ldots,b_r,h)=w_{ij}(-b_1, \ldots,-b_r,h)$. There
are at least two ways to ensure algorithmic identifiability. In a
first approach, given $h$, one may find $(b_1,\ldots,b_r)$ by
minimizing \eqref{4.1}, subject to the constraints $\sum_{k=1}^r
b^2_k=1$ and $b_1 > 0$ (or $b_k>0$ for some $b_k\neq0$ if $b_1=0$). A
second option is to find
$(b_1,\ldots,b_r)$ that minimizes \eqref{4.1} near a suitable
starting point $(c_1,\ldots,c_r)$, satisfying $\sum_{k=1}^r
c^2_k=1$ and $c_1> 0$, and then to rescale the solution to
$ (\frac{b_1}{\sqrt{\sum_k
b^2_k}},\break\ldots,\frac{b_r}{\sqrt{\sum_k
b^2_k}},\frac{h}{\sqrt{\sum_k b^2_k}} )$. The second option is
simpler since the unconstrained minimization is easier to achieve.
However, if one wishes to specify $h$, the constraint
$\sum_{k=1}^r b^2_k=1$ needs to be enforced in the minimization
step. In the simulations, we found that both options led to
virtually the same solution for a well-chosen bandwidth $h$.

The minimization step is a nonlinear least squares problem, which
can be implemented through the optimization package in MATLAB. It
is important to secure a good starting point for the minimization.
We obtained a default starting point by searching along each
dimension separately. Starting with the first dimension, we
located a minimum along $S(b_1)$, as defined in \eqref{4.1}, along
a grid of values of $b_1$ in the interval $[0,1]$. After obtaining
the minimizer $x_1$, we continued to search along the second
dimension using values $S(x_1,b_2)$, where $b_2$ varies on a grid
within $[-1,1]$. This approach was then iterated as necessary and
provided the starting point.

\subsubsection*{Step 2. Selecting $r$ and $h$}

Here $r$ is the number of eigenfunctions used in \eqref{2.8} and
$h$ is the kernel bandwidth. We employed 10-fold cross-validation
to evaluate each pair $(h,r)$. Each of 10 subgroups of curves
denoted by $V_1,\ldots,V_{10}$ was used as a validation set, one
at a time, while the remaining data were used as the training set.
For given $(h, r)$, we found $\hat{\beta}$ as described in step 1 and
computed $S(r,h) =\frac{1}{\sum_k \#V_k}\sum_{k=1}^{10}S_k$, where
$S_k(r,h)=\sum_{j\in V_k}(Y_j-\hat{Y}_j)^2$ and $\hat{Y}_j=\hg
_1(\int
\hat{\be}X_j)$, using local-constant or local-linear method as
described in the paragraph containing \eqref{2.9} and assuming only
$Y_i$ in the training set are known. We then found the
minimizers of $S(r,h)$, which were the selected values
for $r$ and $h$.

\subsubsection*{Step 3. Backfitting step}

By default, we fitted a single-index functional regression model,
which meant that predictions
$\hat{g} (\int\hat{\beta} x_i )$ were obtained via
\eqref{2.5} using the optimal $(h,r)$ chosen in step $2$ and the
corresponding estimated $\beta$ in step $1$. For fitting a $p$-index functional
regression model, the fits obtained in an initial single-index
step gave only $\hat{g}^0_1 (\int\hat{\beta}^0_1 x_i )$ in
\eqref{2.13}. We then replaced $Y_i$ by
$Y_i-\hat{g}^0_1 (\int\hat{\beta}^0_1 x_i )$ and repeated
steps $1$ and $2$ to find
$\hat{g}^0_2 (\int\hat{\beta}^0_2 x_i )$. This procedure
was iterated until $p$ indices were obtained. This only gives us
the initial estimate of the $p$-index model. Then for the $d$th
iteration and the increasing sequence $k=1,\ldots,p$, we used
$Y^d_i(k)$ defined in \eqref{2.14} to fit
$\hg^d_k (\int\hbe^d_k x )$. The iteration stops once
$\|\hbe_1^{d-1}-\hbe_1^{d}\|_{L_2}<0.01$ or $10$ iterations are
reached.\vadjust{\goodbreak}

\subsection{Illustration for spectrometric data}\label{sec4.2}

\begin{figure}

\includegraphics{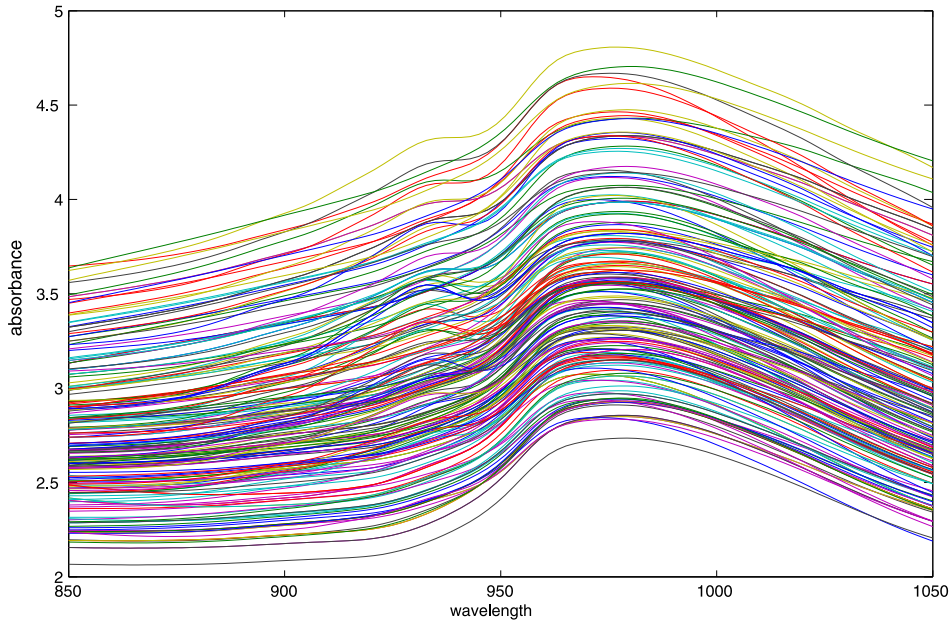}

\caption{Sample of 204 absorbance spectra for meat
specimens.}
\label{fig1}
\end{figure}

We applied the proposed model to spectrometric data that can be
found at \texttt{
\href{http://lib.stat.cmu.edu/datasets/tecator}{http://lib.stat.cmu.edu/}
\href{http://lib.stat.cmu.edu/datasets/tecator}{datasets/tecator}}.
We used only part of the data with data selection performed in the same way
as in~\cite{16} and~\cite{9}. These data were obtained for $215$
pieces of meat, for each of which one observes a spectrometric curve
$X_i$, corresponding to an absorbance spectrum measured at 100
wavelengths. These spectrometric curves are depicted in Figure~\ref{fig1}.
The fat content of each sample was determined by an analytic method
and recorded as a scalar response $Y_i$. One is interested in
predicting the fat content of each sample directly from the
spectrometric curve.

%

%
\begin{figure}

\includegraphics{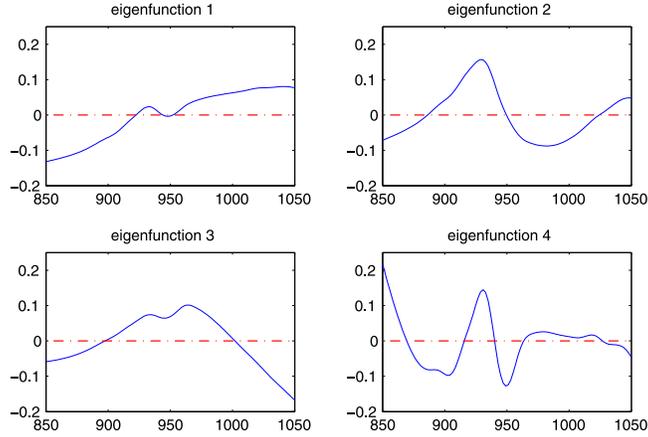}

\caption{The first four estimated eigenfunctions of
the normalized absorbance spectra.}
\label{fig2}
\end{figure}

%
\begin{figure}[b]

\includegraphics{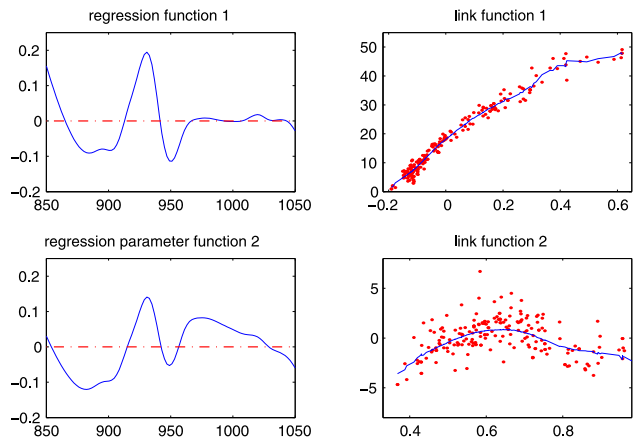}

\caption{The estimated regression parameter
functions and link functions. Left two panels: the estimated
regression parameter functions $\hat{\beta}_1$ and $\hat{\beta}_2$
for the first and second index, respectively; right two panels: the
estimated link functions $\hat{g}_1$ and $\hat{g}_2$ for the first
and second index, respectively.}
\label{fig3}
\end{figure}

In a preprocessing step, we removed $11$ outliers. We also
normalized each spectrometric curve by subtracting its area under
the curve, $\int X_i(t)\,dt$, because we found that the first
eigenfunction of the spectral curves is almost flat and its
eigenvalue is much larger than the others, but the corresponding
fitted coefficient $\hga_1$ in \eqref{2.8} is close to $0$. This
normalization step reduced the leave-one-curve-out prediction error
by more than $30\%$. The first four estimated eigenfunctions for the
normalized curves are plotted in Figure~\ref{fig2}.

To fit the functional single-index model, we used 10-fold
cross-validation to choose the number $r$ of included eigenfunctions in the
representation (\ref{2.8}) and the bandwidth
for the Epanechnikov kernel, obtaining $4$ and $0.0687$ for these
choices. Using the local-linear method described in \eqref{2.5} and
\eqref{2.7}, we then estimated the regression parameter function
$\beta_1$ and the link function~$g_1$. These function estimates are
shown in the upper panel of Figure~\ref{fig3}. The average
leave-one-curve-out squared prediction error for the proposed
single-index model is $3.51$, while fitting a Generalized
Functional Linear Model (GFLM) led to a prediction error of $4.99$,
showing substantial improvement for the proposed model.

We further applied the backfitting procedure described in Section~\ref{sec2.4} to check whether a multiple index functional model is more
appropriate for these data than a single-index model. The average
leave-one-curve-out squared prediction errors were found to be
$2.39$ for the model with two indices and $2.42$ for three indices.
The estimated regression parameter functions $\hat{\beta}_2$ and
link function $\hat{g}_2$ are also displayed in Figure~\ref{fig3}.
The plot
of $\hat{\beta}_1$ suggests that the small bump around wavelength
$930$ is an important indicator of the fat content level. We note
that $\hat{\beta}_2$ has similar shape as $\hat{\beta}_1$ except
for differences around wavelength $975$, where it is positive.
The model with two indices emerges as the best choice for
prediction and improves more than 50\% upon the GFLM and more
than 30\% upon the single-index model in terms of prediction error.

\section{Simulation study}\label{sec5}

\subsection{Simulations for single-index models}\label{sec5.1}

We studied the finite sample performance of five single-index models
(\ref{2.3}). Samples of balanced functional data consisting of
$N=50/200/800$ predictor trajectories and a scalar response were
generated and each predictor function was sampled through 50
equidistantly spaced measurements in $[0, 1]$. The predictor
functions were generated as
\[
X_i(t)=\mu(t)+\sum_{k=1}^4 \xi_{ik}\phi_k(t),  \qquad    i=1,\ldots,N,
\]
where $\mu(t)=t$, $\phi_1(t)=\frac{1}{\sqrt{2}}\sin(2\pi t)$,
$\phi_2(t)=\frac{1}{\sqrt{2}}\cos(2\pi t)$,
$\phi_3(t)=\frac{1}{\sqrt{2}}\sin(4\pi t)$,
$\phi_4(t)=\frac{1}{\sqrt{2}}\cos(4\pi t)$, and $\xi_{ik}$ are
i.i.d. $N(0,\lambda_k)$ with $\lambda_1=1$, $\lambda_2=\frac{1}{2}$,
$\lambda_3=\frac{1}{4}$, $\lambda_4=\frac{1}{8}$. Responses $Y_i$
were obtained as:
\begin{longlist}
\item[Model (i):] $Y_i=\cos (\int_0^1 \beta X_i )+\epsilon_i$
(nonmonotone link);
\item[Model (ii):] $Y_i= (\int_0^1 \beta X_i )^2+\epsilon_i$
(nonmonotone link);
\item[Model (iii):] $Y_i=\int_0^1 \beta X_i+\epsilon_i$ (functional
linear model; trivially, a  monotone link);
\item[Model (iv):] $Y_i \sim\operatorname{Poisson} \{\exp
(2+\int_0^1
\beta X_i ) \}$ (functional generalized
Poisson model; a monotone link with heteroscedastic noise);
\item[Model (v):] $Y_i \sim\operatorname{Binomial} (1,\frac{1}{2}
\cos (2\int_0^1 \beta X_i )+\frac{1}{2} )$
(functional generalized Binomial model; a nonmonotone link with
heteroscedastic noise);
\end{longlist}
where $\beta=\frac{1}{\sqrt{3}}\phi_1+\frac{1}{\sqrt{3}}\phi
_2+\frac{1}
{\sqrt{6}}\phi_3+\frac{1}{\sqrt{6}}\phi_4$
in all models. In models (i)--(iii), errors $\epsilon_i$ were simulated
as i.i.d. Gaussian noise with mean $0$ and $\Var(\epsilon)=R
\Var \{g (\int\beta X ) \}$. Here $R$ is a measure of
the signal-to-noise ratio, with values chosen as $R=0.1$ and $R=0.5$.

We compared the proposed model with the generalized functional
linear regression model (GFLM) with unknown link and variance
function~\cite{13}, which is a single-index model. In the
simulations, we implemented the proposed model using the
local-constant method defined in (\ref{2.4}) (details can be found
in Section~\ref{sec4.1}). Prediction outcomes were quantified by root average
squared errors
$\mathrm{RASE}=\{\frac{1}{N}\sum_i\{\hat{Y}_i-g(\int\beta
X_i)\}^2\}^{1/2}$, where $\hat{Y}_i$ is our estimate of
$g(\int\beta  X_i)$ defined in the paragraph containing \eqref{2.5},
plugging in $\hat{\be}$ and always
leaving $Y_i$ out of the sample when calculating $\hat{Y}_i$. We also
quantified the
error of the estimated regression parameter function by root squared error
$\mathrm{RSE}(\hat{\beta})=\{\int(\hat{\beta}-\beta)^2\}^{1/2}$. Average
values of
RASE and RSE obtained from 100 Monte Carlo runs were then used to evaluate
the procedures.

%
\begin{table}
\tabcolsep=0pt
\caption{Simulation results for single-index models
(\textup{i})--(\textup{iii}). ``FSIR'' denotes the~proposed~functional single-index
regression and ``GFLM'' denotes the~generalized~functional linear
model~\cite{13}}\label{table1}
\begin{tabular*}{\textwidth}{@{\extracolsep{\fill}}lccccccccc@{}}
\hline
   &    & \multicolumn{4}{c}{\textbf{FSIR}}& \multicolumn{4}{c@{}}{\textbf{GFLM}} \\[-5pt]
   &    & \multicolumn{4}{c}{\hrulefill}& \multicolumn{4}{c@{}}{\hrulefill}\\
  & & \multicolumn{2}{c}{$\bolds{R=0.1}$} & \multicolumn{2}{c}{$\bolds{R=0.5}$} & \multicolumn{2}{c}{$\bolds{R=0.1}$} & \multicolumn{2}{c@{}}{$\bolds{R=0.5}$} \\[-5pt]
 & & \multicolumn{2}{c}{\hrulefill} & \multicolumn{2}{c}{\hrulefill} & \multicolumn{2}{c}{\hrulefill} &
 \multicolumn{2}{c@{}}{\hrulefill}\\
\textbf{Model}  &\textit{\textbf{N}} & \textbf{RASE} & \textbf{RSE} & \textbf{RASE} & \textbf{RSE} &
\textbf{RASE} & \textbf{RSE} & \textbf{RASE} &
\textbf{RSE} \\
\hline
 (i)  & \hphantom{8}50 & $ 0.0464$ & $ 0.0204$ & $ 0.1096$ & $ 0.1225$ &
$ 0.1299$ & $ 0.1488$ & $ 0.1546$ & $ 0.2335$ \\
& 200 & $ 0.0279$ & $ 0.0052$ & $ 0.0557$ & $ 0.0195$ & $ 0.0442$ & $ 0.0109$ &
$ 0.0709$ & $ 0.0818$ \\
& 800 & $ 0.0156$ & $ 0.0024$ & $ 0.0315$ & $ 0.0041$ & $ 0.0288$ & $ 0.0025$ &
$ 0.0402$ & $ 0.0049$ \\
 [3pt]
 (ii)  & \hphantom{8}50 & $ 0.1334$ & $ 0.0304$ & $ 0.3071$ & $ 0.2240$ &
$ 0.1914$ & $ 0.1423$ & $ 0.3329$ & $ 0.3407$ \\
& 200 & $ 0.0731$ & $ 0.0065$ & $ 0.1549$ & $ 0.0223$ & $ 0.1058$ & $ 0.0150$ &
$ 0.1838$ & $ 0.0840$ \\
& 800 & $ 0.0399$ & $ 0.0025$ & $ 0.0844$ & $ 0.0047$ & $ 0.0702$ & $ 0.0028$ &
$ 0.0970$ & $ 0.0053$ \\
 [3pt]
 (iii) & \hphantom{8}50 & $ 0.0970$ & $ 0.0341$ & $ 0.2562$ & $ 0.1705$ &
$ 0.1024$ & $ 0.0546$ & $ 0.2378$ & $ 0.1819$ \\
& 200 & $ 0.0486$ & $ 0.0078$ & $ 0.1122$ & $ 0.0332$ & $ 0.0463$ & $ 0.0068$ &
$ 0.1030$ & $ 0.0204$ \\
& 800 & $ 0.0226$ & $ 0.0030$ & $ 0.0526$ & $ 0.0083$ & $ 0.0237$ & $ 0.0026$ &
$ 0.0558$ & $ 0.0071$ \\
\hline
\end{tabular*}
\end{table}

%
\begin{table}[b]
\caption{Simulation results for single-index models (\textup{iv})
and (\textup{v})}\label{table2}
\begin{tabular}{@{}lccccc@{}}
\hline
   &   & \multicolumn{2}{c}{\textbf{FSIR}}& \multicolumn{2}{c@{}}{\textbf{GFLM}} \\[-5pt]
   &   & \multicolumn{2}{c}{\hrulefill}& \multicolumn{2}{c@{}}{\hrulefill}\\
\textbf{Model}  &\textit{\textbf{N}} & \textbf{RASE} & \textbf{RSE} & \textbf{RASE}
& \textbf{RSE} \\
\hline
(iv) & \hphantom{8}50 & $1.798$\hphantom{0} & $ 0.0767$ & $1.632$\hphantom{0} & $ 0.0639$ \\
& 200 & $1.207$\hphantom{0} & $ 0.0214$ & $1.064$\hphantom{0} & $ 0.0137$ \\
& 800 & $ 0.8117$ & $ 0.0071$ & $ 0.6880$ & $ 0.0045$ \\
[3pt]
(v) & \hphantom{8}50 & $ 0.2324$ & $ 0.4023$ & $ 0.2060$ & $ 0.4333$ \\
& 200 & $ 0.1222$ & $ 0.0850$ & $ 0.1400$ & $ 0.2866$ \\
& 800 & $ 0.0612$ & $ 0.0140$ & $ 0.0629$ & $ 0.0728$ \\
\hline
\end{tabular}
\end{table}

The results in Tables~\ref{table1} and~\ref{table2} indicate that the
proposed method works
clearly better than GFLM for models (i), (ii)
and (v), where the link function is nonmonotone. For model (iii),
the performance of the two methods was found to be similar. In this
example, the effect of the monotone link function (here it is
linear) would have been expected to favor the GFLM, but
this may be counteracted by the fact that the GFLM
fits an unnecessarily complex
model in the case of homogeneous errors, as it also includes a nonparametric
variance function estimation step. In model (iv), where the link is
monotone and the noise is heteroscedastic, the GFLM
not unexpectedly
performs better, as it is able to target
the heteroscedastic errors, improving efficiency of the estimates. Overall,
it emerges that the proposed method is
clearly preferable in situations where the link function is
nonmonotone.

\subsection{Simulations for multiple index models}\label{sec5.2}

We simulated data for five multiple index models, using the same
processes and settings as described in Section~\ref{sec5.1}. Three of the
models [(vi)--(viii)] contain two indices and two models [(ix)--(x)]
contain three indices, as follows:
\begin{longlist}
\item[Model (vi):] $Y_i=\cos (\int_0^1 \beta_1 X_i )+0.5\sin
(\int_0^1
\beta_2 X_i )+\epsilon_i$ (two nonmonotone link functions),
where
$\beta_1=\frac{1}{\sqrt{3}}\phi_1+\frac{1}{\sqrt{3}}\phi_2+\frac
{1}{\sqrt{6}}\phi_3+\frac{1}{\sqrt{6}}\phi_4$,
and
$\beta_2=\frac{1}{\sqrt{3}}\phi_1-\frac{1}{\sqrt{3}}\phi_2-\frac
{1}{\sqrt{6}}\phi_3+\frac{1}{\sqrt{6}}\phi_4$;
\item[Model (vii):] $Y_i=\int_0^1 \beta_1 X_i+\exp (0.5\int_0^1
\beta_2 X_i )+\epsilon_i$ (two monotone link functions),
where
$\beta_1=\frac{1}{\sqrt{3}}\phi_1+\frac{1}{\sqrt{3}}\phi_2+\frac
{1}{\sqrt{6}}\phi_3+\frac{1}{\sqrt{6}}\phi_4$
and
$\beta_2=\frac{1}{\sqrt{3}}\phi_1-\frac{1}{\sqrt{3}}\phi_2-\frac
{1}{\sqrt{6}}\phi_3+\frac{1}{\sqrt{6}}\phi_4$;
\item[Model (viii):] $Y_i=\int_0^1 \beta_1 X_i+0.5 (\int_0^1
\beta_2 X_i )^2+\epsilon_i$ (one nonmonotone link and one
monotone link),
where
$\beta_1=\frac{1}{\sqrt{3}}\phi_1+\frac{1}{\sqrt{3}}\phi_2+\frac
{1}{\sqrt{6}}\phi_3+\frac{1}{\sqrt{6}}\phi_4$
and
$\beta_2=\frac{1}{\sqrt{3}}\phi_1-\frac{1}{\sqrt{3}}\phi_2-\frac
{1}{\sqrt{6}}\phi_3+\frac{1}{\sqrt{6}}\phi_4$;
\item[Model (ix):] $Y_i=\int_0^1 \beta_1 X_i+\exp (0.5\int_0^1
\beta_2 X_i )+0.5 (\int_0^1 \beta_1 X_i )^2+\epsilon_i$
(three link functions),
where
$\beta_1=\frac{1}{\sqrt{3}}\phi_1+\frac{1}{\sqrt{3}}\phi_2+\frac
{1}{\sqrt{6}}\phi_3+\frac{1}{\sqrt{6}}\phi_4$,
$\beta_2=\frac{1}{\sqrt{3}}\phi_1-\frac{1}{\sqrt{3}}\phi_2-\frac
{1}{\sqrt{6}}\phi_3+\frac{1}{\sqrt{6}}\phi_4$
and
$\beta_3=-\frac{1}{\sqrt{3}}\phi_1+\frac{1}{\sqrt{3}}\phi
_2+\frac
{1}{\sqrt{6}}\phi_3+\frac{1}{\sqrt{6}}\phi_4$;
\item[Model (x):] $Y_i=\int_0^1 \beta_1 X_i+0.5 (\int_0^1 \beta
_1
X_i )^2+0.25 (\int_0^1
\beta_1 X_i )^3+\epsilon_i$ (three link functions),
where
$\beta_1=\frac{1}{\sqrt{3}}\phi_1+\frac{1}{\sqrt{3}}\phi_2+\frac
{1}{\sqrt{6}}\phi_3+\frac{1}{\sqrt{6}}\phi_4$,
$\beta_2=\frac{1}{\sqrt{3}}\phi_1-\frac{1}{\sqrt{3}}\phi_2-\frac
{1}{\sqrt{6}}\phi_3+\frac{1}{\sqrt{6}}\phi_4$
and
$\beta_3=-\frac{1}{\sqrt{3}}\phi_1+\frac{1}{\sqrt{3}}\phi
_2+\frac
{1}{\sqrt{6}}\phi_3+\frac{1}{\sqrt{6}}\phi_4$.
\end{longlist}

We compared the results from the recursive fitting procedure and the
backfitting procedure in terms of root average squared errors
\[
\mathrm{RASE}_k= \Biggl\{\frac{1}{N}\sum_i \Biggl\{\sum_{j=1}^k\hat
{g}_j
\biggl(\int\hat{\beta}_j
X_i \biggr)-\sum_{j=1}^p g_j \biggl(\int\beta_j
X_i \biggr) \Biggr\}^2 \Biggr\}^{{1}/{2}}
\]
for a $p$-index model
when fitting the first $k$ indices. It is of interest to include
cases $k<p$ (not fitting a sufficient number of indices) and $k>p$
(overfitting the number of indices) and to determine whether the
best results are obtained for the correct number of indices, which
would suggest choosing the number of indices by fitting various
numbers of indices and choosing the number according to the model
with the best fit. Here $\hat{g}_j$ and $\hat{\beta}_j$ are
estimated using both recursive and backfitting procedures. Accordingly,
if the underlying model, selected from models (vi)--(x), contains $p$ indices,
we calculated the values for $\mathrm{RASE}_k$ for $k=1,\ldots,p+1$.

%
\begin{table}
\tabcolsep=0pt
\caption{Simulation results for multiple index
model (\textup{vi})--(\textup{viii}) with two underlying indices.
RASE$^R_k$, $k=1,2,3$, stands for root average errors using the
recursive fitting procedure and $k$ indexes and RASE$^I_k$ for the
same errors obtained when using the iterative backfitting procedure.
Shown are average results based on 100 Monte Carlo runs}\label{table3}
\begin{tabular*}{\textwidth}{@{\extracolsep{\fill}}lcccccccc@{}}
\hline
\textbf{Model}& \textit{\textbf{N}} & \textit{\textbf{R}} & \textbf{RASE}$^{\bolds{R}}_{\bolds{1}}$ & \textbf{RASE}$^{\bolds{R}}_{\bolds{2}}$
& \textbf{RASE}$^{\bolds{R}}_{\bolds{3}}$ & \textbf{RASE}$^{\bolds{I}}_{\bolds{1}}$ &
\textbf{RASE}$^{\bolds{I}}_{\bolds{2}}$ &
\textbf{RASE}$^{\bolds{I}}_{\bolds{3}}$ \\
\hline
(vi) & \hphantom{2}50 & 0.1 & $ 0.2975$ & $ 0.1960$ & $ 0.1872$ & $ 0.2975$
& $ 0.1209$ & $ 0.1483$ \\
& 200 & 0.1 & $ 0.3003$ & $ 0.1357$ & $ 0.1059$ & $ 0.3003$ & $ 0.0645$ &
$ 0.0854$ \\
& \hphantom{2}50 & 0.5 & $ 0.3206$ & $ 0.2683$ & $ 0.2797$ & $ 0.3206$ & $ 0.2048$ & $ 0.2810$
\\
& 200 & 0.5 & $ 0.3051$ & $ 0.1778$ & $ 0.1754$ & $ 0.3051$ & $ 0.1235$ &
$ 0.1427$ \\
[3pt]
(vii) & \hphantom{2}50 & 0.1 & $ 0.2107$ & $ 0.2076$ & $ 0.2144$ & $ 0.2107$
& $ 0.1991$ & $ 0.2312$ \\
& 200 & 0.1 & $ 0.1311$ & $ 0.1131$ & $ 0.1372$ & $ 0.1311$ & $ 0.1070$ &
$ 0.1247$ \\
& \hphantom{2}50 & 0.5 & $ 0.4238$ & $ 0.3937$ & $ 0.4592$ & $ 0.4238$ & $ 0.3786$ & $ 0.4335$
\\
& 200 & 0.5 & $ 0.2507$ & $ 0.2329$ & $ 0.2719$ & $ 0.2507$ & $ 0.2267$ & $ 0.2848$
\\
 [3pt]
(viii) & \hphantom{2}50 & 0.1 & $ 0.4463$ & $ 0.3184$ & $ 0.3317$ &
$ 0.4463$ & $ 0.2310$ & $ 0.2817$ \\
& 200 & 0.1 & $ 0.4061$ & $ 0.1496$ & $ 0.1594$ & $ 0.4061$ & $ 0.1016$ &
$ 0.1279$ \\
& \hphantom{2}50 & 0.5 & $ 0.4818$ & $ 0.4872$ & $ 0.5327$ & $ 0.4818$ & $ 0.4632$ & $ 0.4698$
\\
& 200 & 0.5 & $ 0.4304$ & $ 0.2502$ & $ 0.2910$ & $ 0.4304$ & $ 0.2107$ &
$ 0.2412$ \\
\hline
\end{tabular*}
\end{table}
%
\begin{table}[b]
\caption{Recursive fitting results for
models (\textup{ix}) and (\textup{x}) with three indices}\label{table4}
\begin{tabular}{@{}lcccccc@{}}
\hline  \textbf{Model} & \textit{\textbf{N}} & \textit{\textbf{R}} & \textbf{RASE}$^{\bolds{R}}_{\bolds{1}}$ & \textbf{RASE}$^{\bolds{R}}_{\bolds{2}}$ & \textbf{RASE}$^{\bolds{R}}_{\bolds{3}}$ &
\textbf{RASE}$^{\bolds{R}}_{\bolds{4}}$ \\
\hline
(ix) & \hphantom{2}50 & 0.1 & $ 0.5107$ & $ 0.3417$ & $ 0.3518$ & $ 0.3772$
\\
& 200 & 0.1 & $ 0.4791$ & $ 0.2196$ & $ 0.2132$ & $ 0.2183$ \\
& \hphantom{2}50 & 0.5 & $ 0.5453$ & $ 0.4810$ & $ 0.5104$ & $ 0.5297$\\
& 200 & 0.5 & $ 0.5161$ & $ 0.3324$ & $ 0.3329$ & $ 0.3504$\\
 [3pt]
(x) & \hphantom{2}50 & 0.1 & $ 0.5107$ & $ 0.3417$ & $ 0.3518$ & $ 0.3772$
\\
& 200 & 0.1 & $ 0.4792$ & $ 0.2327$ & $ 0.2264$ & $ 0.2316$ \\
& \hphantom{2}50 & 0.5 & $ 0.6631$ & $ 0.6461$ & $ 0.6111$ & $ 0.6418$ \\
& 200 & 0.5 & $ 0.5161$ & $ 0.3224$ & $ 0.3329$ & $ 0.3504$ \\
\hline
\end{tabular}
\end{table}

%
\begin{table}
\tabcolsep=0pt
\caption{Iterative backfitting results for
models (\textup{ix}) and (\textup{x}) with three indices}\label{table5}
\begin{tabular*}{260pt}{@{\extracolsep{\fill}}lcccccc@{}}
\hline  \textbf{Model} & \textit{\textbf{N}} & \textit{\textbf{R}} & \textbf{RASE}$^{\bolds{I}}_{\bolds{1}}$
& \textbf{RASE}$^{\bolds{I}}_{\bolds{2}}$ & \textbf{RASE}$^{\bolds{I}}_{\bolds{3}}$ &
\textbf{RASE}$^{\bolds{I}}_{\bolds{4}}$ \\
\hline
(ix) & \hphantom{2}50 & 0.1 & $ 0.5107$ & $ 0.3272$ & $ 0.3009$ & $ 0.3395$
\\
& 200 & 0.1 & $ 0.4791$ & $ 0.2084$ & $ 0.1422$ & $ 0.1988$ \\
& \hphantom{2}50 & 0.5 & $ 0.5453$ & $ 0.5372$ & $ 0.5518$ & $ 0.6095$ \\
& 200 & 0.5 & $ 0.5161$ & $ 0.3350$ & $ 0.3137$ & $ 0.3980$ \\
 [3pt]
(x) & \hphantom{2}50 & 0.1 & $ 0.5107$ & $ 0.3501$ & $ 0.3106$ & $ 0.3495$
\\
& 200 & 0.1 & $ 0.4792$ & $ 0.2063$ & $ 0.1808$ & $ 0.1860$ \\
& \hphantom{2}50 & 0.5 & $ 0.6631$ & $ 0.6291$ & $ 0.5825$ & $ 0.6476$ \\
& 200 & 0.5 & $ 0.5161$ & $ 0.3372$ & $ 0.3129$ & $ 0.3248$ \\
\hline
\end{tabular*}
\end{table}

As one can see from the results in Tables~\ref{table3},~\ref{table4}
and~\ref{table5}, the recursive fitting procedure
often does not identify the right number of indexes and for nearly
all fits produces larger RASE values, as compared to the iterative
backfitting procedure. The iterative
backfitting
 method thus emerges as the
preferred method.

\begin{appendix}\label{appm}

\section*{\texorpdfstring{Appendix: Proof of Theorem~\lowercase{\protect\ref{theo3.1}}}{Appendix: Proof of Theorem 3.1}}
\label{app}

We describe the details of the proof by breaking it up into
several steps.

\subsection*{Step~1. Upper bound on mean summed squared error}

Define $\ga_j=g(X_j)=g_1(\inti\be\z X_j)$ and
%
\begin{equation}\label{A.1} \qquad
\bga_j= \biggl(\suminj\ga_i K\ij \biggr) \Big/\suminj K\ij , \qquad
\bep_j= \biggl(\suminj\ep_i K\ij \biggr)\Big /\suminj K\ij .\
\end{equation}
To
express their dependence on $\be$, through $K\ij=K\ij(\be)$, we
shall write $\bga_j$, $\bep_j$ and $\bY_j$ as $\bga(\be)$,
$\bep_j(\be)$ and $\bY_j(\be)$, respectively. In this notation,
$S(\be)=S_0+S_1(\be)+S_2(\be)+2 S_3(\be)$, where $S_0=\sum_{1\leq
j\leq n} \ep_j^2$ and does not depend on $\be$,
%
\begin{eqnarray}\label{A.2}
 S_1(\be)&=&\sumjon\{\ga_j-\bY_j(\be)\}^2 ,\qquad
S_2(\be)=\sumjon\{\ga_j-\bga_j(\be)\} \ep_j ,\nonumber
\\[-8pt]
\\[-8pt]
 S_3(\be)&=&\sumjon\bep_j(\be) \ep_j .
\nonumber
\end{eqnarray}
Furthermore,
$S_1(\be)=S_4(\be)-2 S_5(\be)+S_6(\be)$, where
%
\begin{eqnarray}\label{A.3}
 S_4(\be)&=&\sumjon\{\ga_j-\bga_j(\be)\}^2 ,\qquad
S_5(\be)=\sumjon\{\ga_j-\bga_j(\be)\} \bep_j(\be) ,\nonumber
\\[-8pt]
\\[-8pt]
 S_6(\be)&=&\sumjon\bep_j(\be)^2 ,
\nonumber
\end{eqnarray}
with notations as in (\ref{A.1}).

Let $\cB_1=\cB_1(n)$ denote a class of functions $\be$, and
suppose we can prove that
%
\begin{equation}
\sup_{\be\in\cB_1} |S_k(\be)|=O_p({\lambda }_n)  \qquad\mbox
{for }
k=2,3,5,6 ,\label{A.4}
\end{equation}
where ${\lambda }_n$ denotes a sequence of
positive constants. Then,
%
\begin{eqnarray}\label{A.5}
S_1(\hbe)
&=&S(\hbe)-\{S_0+S_2(\hbe)+2 S_3(\hbe)\}\nonumber\\
&\leq&
S(\be\z)-\{S_0+S_2(\hbe)+2 S_3(\hbe)\}\nonumber\\
&=&
S_1(\be\z)+S_2(\be\z)+2 S_3(\be\z) -\{S_2(\hbe)+2 S_3(\hbe
)\}\nonumber
\\[-8pt]
\\[-8pt]
&=&
S_4(\be\z)-2 S_5(\be\z)+S_6(\be\z) +S_2(\be\z)
\nonumber\\
&&{}         +2 S_3(\be
\z)-\{
S_2(\hbe)+2 S_3(\hbe)\}
\nonumber\\&=&
S_4(\be\z)+O_p({\lambda }_n) ,
\nonumber
\end{eqnarray}
where the inequality follows from the fact that $\be=\hbe$ minimizes
$S(\be)$, the final identity follows from \eqref{A.4} provided that
$\be\z$ and $\hbe$ are both in $\cB_1(n)$, and all other identities
in this string hold true generally.

Without loss of generality, the support of $K$ is contained in the
interval $[-1,1]$ [see \eqref{3.7}]. If in addition
$|g_1(u)-g_1(v)|\leq D_1 |u-v|^{a_1}$ for all $u$ and $v$ [see
\eqref{3.2}], then $|\ga_j-\bga_j(\be\z)|\leq D_1 h^{a_1}$ for all
$j$, and therefore
%
\begin{equation}S_4(\be\z)\leq
n  (D_1 h^{a_1} )^2 .\label{A.6}
\end{equation}
Together, \eqref{A.5}
and \eqref{A.6} imply that
%
\begin{equation}\sumjon\{g(X_j)-\hg(X_j)\}^2
=O_p ({\lambda }_n+n h^{2a_1} ) .\label{A.7}
\end{equation}

\subsection*{\texorpdfstring{Step~2. Decomposition of each set $S_k(\be)$ into two parts}
{Step~2. Decomposition of each set $S_k(beta)$ into two parts}}

Let $\cX=\{X_1,\ldots,\break X_n\}$ denote the set of explanatory
variables, and for each $\be\in\cB_1$ let
$\cJ=\cJ(\be)\subseteq\cJ\z\equiv\{1,\ldots,n\}$ denote a
random set
which satisfies
%
\begin{equation}
P [\#\{\cJ\z\setminus\cJ(\be)\}>2 D_5 n h^{a_4}  \mbox{ for
some } \be\in\cB_1 ]\ra0\label{A.8}
\end{equation}
as $n\rai$, where $a_4$
is as in (\ref{3.6}). (The set $\cJ$ will be $\cX$-measurable.) Define
$S_k\J(\be)$, for $2\leq k\leq6$, to be the version of $S_k(\be)$
that arises if, in the definitions at \eqref{2.1} and \eqref{2.2},
we replace summation over $1\leq j\leq n$ by summation over
$j\in\cJ$. Since $g$ is bounded, and all moments of the error
variables $\ep_i$ are finite [see \eqref{3.3}], then $\sup_{1\leq
i\leq n} |Y_i|=O_p(n^\eta)$ with probability~1, for all $\eta>0$.
Therefore, in view of~\eqref{A.8},
%
\begin{equation}\quad
\max_{k=1,\ldots,6} \sup_{\be\in\cB_1}  |S_k(\be)-S_k\J
(\be) |
=O_p (n^{1+\eta} h^{a_4} )\qquad\mbox{for all }
\eta>0 .\label{A.9}
\end{equation}

\subsection*{\texorpdfstring{Step~3. Determining $\cJ$ for which (\protect\ref{A.8}) holds}
{Step~3. Determining $\cJ$ for which (A.8) holds}}

Define $T_j(\be)=\suminj K\ij$, recall that $f( \cdot\mi\be)$
denotes the probability density of $\inti\be X$, and put
\[
\a_j(\be)=h\int K(u) f \biggl({\inti\be
X_j}-hu  \mid \be \biggr)\,du .
\]
Then,
\begin{eqnarray*}
E\{T_j(\be)\mi X_j\} &=&\a_j(\be) ,\\
\operatorname{var}\{T_j(\be
)\mi X_j\}
&\leq&(n-1) h\int K(u)^2 f\biggl ({\inti\be
X_j}-hu  \mid \be \biggr)\,du\\ &\leq& n (\sup K) \a_j(\be) .
\end{eqnarray*}
Moreover, $0\leq K\ij\leq\sup K$. Therefore by Bernstein's
inequality, if \mbox{$0<c_1<1$},
%
\begin{eqnarray}\label{A.10}
&&P\{T_j(\be)\leq(1-c_1) n \a_j(\be)\mi X_j\}\nonumber\\
&& \qquad =P\{n \a_j(\be)-T_j(\be)\geq c_1 n \a_j(\be)\mi X_j\}\nonumber
\\[-8pt]
\\[-8pt]
&& \qquad \leq\exp\biggl [-\frac{\{c_1 n \a_j(\be)\}^2/2 }{(\sup
K) \{n \a_j(\be)+\otd c_1 n \a_j(\be)\}} \biggr]\nonumber\\
&& \qquad =\exp \biggl\{-\frac{c_1^2 n \a_j(\be)}{2 (\sup
K) (1+\otd c_1)} \biggr\} .
\nonumber
\end{eqnarray}
Hence, defining $\cJ(\be)$ to be the set of all integers $j$ such
that $\a_j(\be)\geq n^{-c_2} h$, where $0<c_2<1$; and putting
$C_2=c_1^2/\{2 (\sup K) (1+\frac13 c_1)\}$; we obtain
\[
\sup_{j\in\cJ(\be)} P\{T_j(\be)\leq(1-c_1) n \a_j(\be)\mi
X_j\}
\leq\exp (-C_2 n^{1-c_2} h ) .
\]
Therefore, since
$\cJ(\be)$ contains no more than $n$ elements, then
\begin{eqnarray*}
&&P \{T_j(\be)\leq(1-c_1) n \a_j(\be) \mbox{ for some }
j\in\cJ(\be)\mbox{ and some } \be\in\cB_1 \}\\
&& \qquad \leq
n (\#\cB_1) \exp (-C_2 n^{1-c_2} h ) .
\end{eqnarray*}
Hence, provided
%
\begin{equation}
\#\cB_1=O \{n^{-C_3-1} \exp (C_2 n^{1-c_2} h ) \}
\label{A.11}
\end{equation}
for some $C_3>0$, we have
%
\begin{equation} \qquad P \{
T_j(\be)>(1-c_1) n \a_j(\be) \mbox{ for all } j\in\cJ(\be)
\mbox{ and all } \be\in\cB_1  \} \ra1 .\label{A.12}
\end{equation}

Note, too, that if $a_3$ and $a_4$ are as in \eqref{3.6}, if $K$ is
supported on $[-1,1]$, and if
%
\begin{equation}(\sup K)\mo n^{-c_2}\leq
D_4 h^{a_3} ,\label{A.13}
\end{equation}
then
\begin{eqnarray*}
\#\{\cJ\z\setminus\cJ(\be)\} &=&\sumjon
I \{\a_j(\be)<n^{-c_2} h \}\nonumber\\
 &=&\sumjon I \biggl\{\inti
K(u)
f \biggl(\inti\be X_j-hu\biggmi\be \biggr)\,du<n^{-c_2} \biggr\}
\\
&\leq&\sumjon I_j ,\nonumber
\end{eqnarray*}
where
\[I_j=I_j(\be)=I \biggl\{\sup_{|u|\leq h}  f
\biggl(\inti\be
X_j-u\biggmi\be \biggr)<D_4 h^{a_3} \biggr\} .
\]
The random
variables $I_1,\ldots,I_n$ are independent and identically
distributed, and, in view of \eqref{3.6}, $\pi(\be)\equiv
P\{I_j(\be)=1\}\leq D_5 h^{a_4}$. Therefore, by Bernstein's
inequality,
\begin{eqnarray*}
P \Biggl\{ \sumjon I_j(\be)>2 D_5 h^{a_4} \Biggr\} &\leq&
P \Biggl[\sumjon\{I_j(\be)-\pi(\be)\}>D_5 h^{a_4} \Biggr]\\
&\leq&\exp \biggl[-\frac{(D_5 h^{a_4})^2/2}{
n \pi(\be) \{1-\pi(\be)\}+\otd D_5 h^{a_4}} \biggr]\\
&\leq&\exp (-3 D_5 n h^{a_4} /8 ) .
\end{eqnarray*}
Hence, provided
%
\begin{equation}
\#\cB_1=o \{\exp (3 D_5 n h^{a_4} /8 ) \}
,\label{A.14}
\end{equation}
result \eqref{A.8} holds.

\subsection*{\texorpdfstring{Step~4. Bound for $E_\cX\{S_k^\cJ(\be)^{2m}\}$ for
$k=2,3,5,6$ and integers \mbox{$m\geq1$}}
{Step~4. Bound for $E_\cX\{S_k^\cJ(beta)^{2m}\}$ for
$k=2,3,5,6$ and integers m >= 1}}

Write $E_\cX$ for expectation conditional on $\cX$, let $Q=Q(\be)$
denote the infimum of $\suminj K\ij$ over all $j\in\cJ$, and put
$\si^2=E(\ep^2)$. Defining $L\ij=K\ij/(\sum_{i_1 \dvtx  i_1\neq
j} K_{i_1j})$, taking $m\geq1$ to be an integer, and using
Rosenthal's inequality, we deduce that for a constant $A(m)$
depending only on $m$,
%
\begin{eqnarray}\label{A.15}
 E_\cX (\bep_j^{2m} ) &\leq&
A(m)  \biggl\{\si^{2m}  \biggl(\sum_{i \dvtx  i\neq j}
L_{ij}^2 \biggr)^{ m} +E (\ep^{2m} ) \sum_{i \dvtx  i\neq j}
L_{ij}^{2m} \biggr\}\nonumber
\\[-4pt]
\\[-12pt]
&\leq& A(m)  \bigl\{ (\si^2 Q\mo \sup
K )^m +E (\ep^{2m} ) Q^{-(2m-1)} (\sup
K)^{2m-1} \bigr\} .
\nonumber
\end{eqnarray}
Therefore,
%
\begin{eqnarray}\label{A.16}
&&E_\cX  \{S_6\J(\be)^{2m} \}\nonumber\\
&& \qquad \leq \biggl\{\sumcj (E_\cX|\bep_j|^{2m} )^{1/(2m)} \biggr\}
^{ 2m}
\\
&& \qquad \leq A(m) n^{2m} \bigl \{ (\si^2 Q\mo \sup K )^m
+E (\ep^{2m} ) Q^{-(2m-1)} (\sup
K)^{2m-1} \bigr\} .
\nonumber
\end{eqnarray}
Moreover, if $|g|\leq C_1$, then $S_4\J(\be)\leq n (2C_1)^2$, and
so, since $S_5\J(\be)^2\leq S_4\J(\be)\times  S_6\J(\be)$, then
%
\begin{eqnarray}\label{A.17}
E_\cX \{S_5\J(\be)^{2m} \} &\leq& S_4\J(\be)^m  \{
E_\cX
S_6\J(\be)^{2m} \}\half\nonumber
\\[-8pt]
\\[-8pt]
&\leq& \{n (2C_1)^2 \}^m
 \{E_\cX S_6\J(\be)^{2m} \}\half .
\nonumber
\end{eqnarray}
More
simply, if $|g|\leq C_1$, then $S_4\J(\be)\leq n (2C_1)^2$ and
$\sumj|\ga_j-\bga_j(\be)|^{2m}\leq n (2C_1)^{2m}$, both uniformly
in $\be$. Therefore,
%
\begin{eqnarray} \label{A.18}
E_\cX \{S_2\J(\be)^{2m} \} &\leq&
A(m)  \biggl[ \{\si^2 S_4\J(\be) \}^m
+E (\ep^{2m} ) \sumcj|\ga_j-\bga_j(\be)|^{2m} \biggr]\nonumber
\\[-8pt]
\\[-8pt]
&\leq&
A(m) (2C_1)^{2m}  \{(n \si^2)^m
+n E (\ep^{2m} ) \} .
\nonumber
\end{eqnarray}

Recall that the support of $K$ is contained in the interval
$[-1,1]$. Let $N_1$ denote the maximum, over values $j\in\cJ$, of
the number of indices $k$ such that $|\inti\be (X_j-X_k)|\leq h$.
Then, the series $\suminj K\ij$ has, for each $j$, at most $N_1$
nonzero terms. Array the values of $\inti\be X_j$, for $j\in\cJ$,
on the real line, and group them into consecutive blocks of indices
$j$, each block (except for the last remnant block) containing just
$N_1$ values. Index these blocks, from left to right along the
line, from 1 to $N_2$, where $N_2$ equals
$\lfloor(\#\cJ)/N_1\rfloor$ or $\lfloor(\#\cJ)/N_1\rfloor+1$ and
$\lfloor x\rfloor$ denotes the integer part of~$x$. Choose one point
$\inti\be X_j$ from each even-indexed block, and remove those
points from the respective blocks; and repeat this until all the
points are removed from all the blocks. Record, for each pass
through the $N_2$ blocks, the removed sequence $j_1,\ldots,j_\nu$ of
indices. (On the first pass, $\nu$ will equal $\lfloor N_2/2\rfloor$
or $\lfloor N_2/2\rfloor+1$, but on later passes, $\nu$ may be
reduced in size.) Now repeat this for odd-indexed blocks. Denote by
$j_{k1},\ldots,j_{kM_k}$, for $1\leq k\leq N$ say, the different
sequences $j_1,\ldots,j_\nu$ that are obtained in this way. The set
of all such sequences represents a (disjoint) partition of the
integers in $\cJ$, and in particular, $M_1+\cdots+M_N=n$. By
construction, for each $k$ the random variables
$\ep_{j_{k1}} \bep_{j_{k1}},\ldots,\ep_{j_{kM_k}} \bep_{j_{kM_k}}$
are independent, conditional on $\cX$; the random integers $N$ and
$M_1,\ldots,M_N$ are measurable in the sigma-field generated by
$\cX$; $N\leq2N_1$; and
$\max_k M_k\leq\lfloor(\#\cJ)/(2N_1)\rfloor+1$.

Since
\[\sumcj\ep_j \bep_j
=\sum_{k=1}^N  (\ep_{j_{k1}} \bep_{j_{k1}}
+\cdots +\ep_{j_{kM_k}} \bep_{j_{kM_k}} ),
\]
then, for any
integer $m\geq1$ and an absolute constant $A(m)\geq1$, depending
only on~$m$,
\begin{eqnarray*}
&&\hspace*{-2pt}E_\cX  \{S_3\J(\be)^{2m} \}\\
&&\hspace*{-2pt} \qquad =E_\cX \biggl\{ \biggl(\sumcj\ep_j \bep_j \biggr)^{ 2m} \biggr\}\\
&&\hspace*{-2pt} \qquad \leq \Biggl(\sum_{k=1}^N  [E_\cX \{ |\ep_{j_{k1}}
\bep
_{j_{k1}}+\cdots
+\ep_{j_{kM_k}} \bep_{j_{kM_k}} |^{2m} \}
]^{1/(2m)} \Biggr)^{
2m}\\
&&\hspace*{-2pt} \qquad \leq A(m)  \Biggl(\sum_{k=1}^N  \Biggl[ \Biggl\{\sum_{\ell=1}^{M_k}
E_\cX (\ep_{j_{k\ell}}^2 \bep_{j_{k\ell}}^2 ) \Biggr\}^{ m}
+\sum_{\ell=1}^{M_k}  E_\cX (|\ep_{j_{k\ell}}
\bep_{j_{k\ell}}|^{2m} ) \Biggr]^{1/(2m)} \Biggr)^{ 2m}\\
&&\hspace*{-2pt} \qquad \leq
A(m) N^{2m}\\
&&\hspace*{-2pt} \qquad  \quad {}\times \max_{1\leq k\leq N}  \Bigl[\si^{2m}
 \Bigl\{M_k \max_{1\leq\ell\leq M_k}
E_\cX (\bep_{j_{k\ell}}^2 ) \Bigr\}^m
+E\ep^{2m} M_k \max_{1\leq
\ell\leq M_k}  E_\cX (|\bep_{j_{k\ell}}|^{2m} ) \Bigr] .
\end{eqnarray*}
Therefore, by (\ref{A.15}),
%
\begin{eqnarray}\label{A.19}
&&E_\cX \{S_3\J(\be)^{2m} \}\nonumber\\ && \qquad \leq
A(m)^2 N^{2m}  \Bigl\{ (\si^2 Q\mo \sup K )^m  \max
_{1\leq
k\leq N} M_k^m\\
&& \qquad  \quad \hphantom{A(m)^2 N^{2m}  \Bigl\{}{}
+E (\ep^{2m} ) Q^{-(2m-1)} (\sup K)^{2m-1}  \max_{1\leq
k\leq N} M_k \Bigr\} .
\nonumber
\end{eqnarray}

The constant $A(m)$ in these bounds can be taken equal to
$(A m/\log m)^m$, where $A>1$ denotes an absolute constant
\cite{24,25}. From this property, and results \eqref{A.16},
\eqref{A.17}, \eqref{A.18} and \eqref{A.19}, and recalling that
$N\leq2 N_1$ and $M_k\leq(n/2N_1)+1$, we deduce that for a
constant $C_4>1$,
%
\begin{eqnarray}\label{A.20}
&&\sum_{k=2,3,5,6}E_\cX  \{S_k\J(\be)^{2m} \}\nonumber\\&& \qquad \leq(m/\log
m)^{m/2} (C_4n)^{2m}  \bigl\{Q^{-m/2}
+E (\ep^{2m} ) Q^{(1/2)-m} \bigr\}\\
&& \qquad  \quad {}+(m/\log
m)^{2m} C_4^m  \{(n N_1/Q)^m
+E (\ep^{2m} ) n (N_1/Q)^{2m-1} \} .
\nonumber
\end{eqnarray}
The contributions to the left-hand side from $S_3\J$ and $S_5\J$
dominate, and so the right-hand side represents, in effect,
$E_\cX\{S_3\J(\be)^{2m}\}+E_\cX\{S_5\J(\be)^{2m}\}$.

\subsection*{Step~5. Upper bounds for $N_1$ and $Q\mo$}

Let $T_j\one(\be)$ denote the version of $T_j(\be)$ in the special
case where $K\equiv1$ on $[-1,1]$ and $K=0$ elsewhere, and write
$\a_j\one(\be)=h\int_{|u|\leq1} f(\inti\be X_j-hu\mi\be)\,du$,
representing the corresponding version of $\a_j(\be)$. In this
notation, $N_1=N_1(\be)$ equals the maximum, over $j$, of the values
of $T_j\one(\be)$ for $j\in\cJ(\be)$. The argument leading to
\eqref{A.10} now gives
\begin{eqnarray*}
&&P \bigl\{T_j\one(\be) >(1+c_1) n \a_j\one(\be)\Bigmi X_j \bigr\}\\
&& \qquad =P \bigl\{T_j\one(\be)-n \a_j\one(\be)\geq
c_1 n \a_j\one(\be)\Bigmi X_j \bigr\}\\
&& \qquad \leq\exp \biggl[-\frac{\{c_1 n \a_j\one(\be)\}^2/2 }{
n \a_j\one(\be)+\otd c_1 n \a_j\one(\be)} \biggr]
=\exp \biggl\{-\frac{c_1^2 n \a_j\one(\be)
}{2 (1+\otd c_1)} \biggr\} .
\end{eqnarray*}
The analogue of \eqref{A.12} in this setting is, assuming that
\eqref{A.11} holds:
%
\begin{equation}P \bigl\{
T_j\one(\be)\leq(1+c_1) n \a_j\one(\be) \mbox{ for all }j
\mbox{ and all } \be\in\cB_1  \bigr\} \ra1 .\label{A.21}
\end{equation}
Since
$\a_j\one(\be)\leq h \sup f(\cdot\mi\be)$, then, using \eqref{3.5},
we deduce from \eqref{A.21} that for a constant $C_5>0$,
%
\begin{equation}
P \{N_1(\be)\leq C_5 nh \mbox{ for all } \be\in\cB_1
\}
\ra1 .\label{A.22}
\end{equation}

Observe, too, that
%
\begin{eqnarray}\label{A.23}
Q(\be)\mo&= &\Bigl\{\inf_{j\in\cJ(\be)} T_j (\be) \Bigr\}\mo
\leq\Bigl \{(1-c_1) \inf_{j\in\cJ(\be)} n \a_j(\be) \Bigr\}\mo
\nonumber
\\[-8pt]
\\[-8pt]
& \leq&
(1-c_1)\mo n^{c_2-1} h\mo ,
\nonumber
\end{eqnarray}
where the first
identity is just the definition of $Q$; the second, in view of
(\ref{A.12}), holds uniformly in $\be\in\cB_1$, with probability
converging to 1 as $n\rai$; and the third is a consequence of the
definition of $\cJ(\be)$ as the set of $j$ for which $\a_j(\be)\geq
n^{-c_2} h$.

\subsection*{\texorpdfstring{Step~6. Proof of uniform convergence to zero of
$n\mo S_k\J(\be)$ for $k=2,3,5,6$}
{Step~6. Proof of uniform convergence to zero of
$n\mo S_k\J(beta)$ for $k=2,3,5,6$}}

Incorporating the bounds at \eqref{A.21} and \eqref{A.22} into
\eqref{A.20}, and taking $m$ to diverge polynomially fast in $n$, we
deduce that, for constants $C_6,C_7>1$, and with probability
converging to 1 as $n\rai$,
%
\begin{eqnarray}\label{A.24}
s(m,n)&\equiv&\sup_{\be\in\cB_1} \sum_{k=2,3,5,6}E_\cX \{
S_k\J(\be
)^{2m} \}\nonumber\\
&\leq&(m/\log
m)^{m/2} (C_6n)^{2m}  \bigl\{ (n^{c_2-1} /h )^{m/2}+E
(\ep
^{2m} )   (n^{c_2-1} /h )^{m-(1/2)} \bigr\}\nonumber\\
&&{}+(m/\log
m)^{2m}
C_6^m  \bigl\{n^{m(c_2+1)}   +E (\ep^{2m} ) n^{(2m-1)c_2+1} \bigr\}\\
&\leq&(C_7 n)^{2m}  \{ (mn^{c_2-1} /h )^{m/2}
+ (m^{2a_2}n^{c_2-1} /h )^{m} \}\nonumber\\
&&{}  + (C_7 n^2 )^m  \bigl\{ (m^2 n^{c_2-1} )^m
+ (m^{2(a_2+1)} n^{2c_2-2} )^m \bigr\} ,
\nonumber
\end{eqnarray}
where, to obtain the last inequality, we used the bound
$E|\ep|^m\leq(D_2 m)^{a_2m}$ in \eqref{3.3}.

Choose $c_2$, and further positive constants $C_8,C_9,c_3,c_4,c_5$,
such that
%
\begin{equation}c_2+2c_3 \max(1,a_2)+c_5<1\quad\mbox{and} \quad
0<c_4<c_5 .\label{A.25}
\end{equation}
Take $m$ equal to the integer part of
$n^{c_3}$ and
%
\begin{equation}C_8 n^{-c_5}\leq h\leq
C_9 n^{-c_4} .\label{A.26}
\end{equation}
The constant $c_2\in(0,1)$ was
introduced immediately below \eqref{A.10}, and, up to \eqref{A.25},
was subject only to the conditions \eqref{A.13} and $0<c_2<1$. For
any given $a_3$ and $c_2$, no matter how small the latter, we can
ensure that \eqref{A.13} holds merely by taking $c_5$ (and thence
$c_4$), in \eqref{A.26}, sufficiently small. Since the results
below continue to hold no matter how small we choose $c_5$ (and
$c_4$), then we can be sure that \eqref{A.13} is satisfied.

It follows from \eqref{A.24}--\eqref{A.26} that, with probability
converging to 1 as ${n\rai}$,
\begin{eqnarray*}
s(m,n)&\leq&(C_7 n)^{2m}  \bigl\{ (n^{c_2+c_3+c_5-1} )^m
+ (n^{c_2+2a_1c_3+c_5-1} )^m\\
&&\hphantom{(C_7 n)^{2m}  \bigl\{}{} + (n^{c_2+2c_3-1} )^m
+ (n^{2\{c_2+c_3(a_2+1)+c_2-1\}} )^m \bigr\}\\
&\leq&4  (C_7 n^{c_6+1} )^m ,
\end{eqnarray*}
where
$c_6=\max\{c_2+c_3+c_5,c_2+2a_1c_3+c_5,c_2+2c_3,c_2+c_3(a_2+1)+ c_2\}<1$.
Therefore, if $0<c_7<1-c_6$ and we put $c_8=(1-c_6-c_7)/2>0$, then,
by Markov's inequality,
%
\begin{eqnarray}\label{A.27}
 \qquad P \biggl\{n\mo \sup_{\be\in\cB_1} \sum_{k=2,3,5,6}|S_k\J(\be)|
>n^{-c_8}\biggmi\cX \biggr\}
&\leq&16^m (\#\cB_1) s(m,n) n^{-2m(1-c_8)}\nonumber
\\[-4pt]
\\[-12pt]
 \qquad &\leq&4 (\#\cB_1)  (16 C_7 n^{-c_7} )^m ,
\nonumber
\end{eqnarray}
where the inequalities hold with probability converging to 1 as
$n\rai$. Hence, provided that
%
\begin{equation}
(\#\cB_1)  (16 C_7 n^{-c_7} )^m\ra0 ,\label{A.28}
\end{equation}
the
left-hand side of \eqref{A.27} converges in probability to zero as
$n\rai$. It follows that the unconditional form of that probability
also converges to zero, and hence that
%
\begin{equation}
P \Bigl\{n\mo \sup_{\be\in\cB_1} \max_{k=2,3,5,6}|S_k\J(\be
)|>n^{-c_8} \Bigr\}
\ra0 .\label{A.29}
\end{equation}

\subsection*{Step~7. Completion}

From \eqref{A.9} and \eqref{A.29} we deduce that, for all $\eta>0$,
%
\begin{equation}
  n\mo \sup_{\be\in\cB_1} \max_{k=2,3,5,6}|S_k(\be)|
=O_p (n^\eta h^{a_4}+n^{-c_8} )
=O_p (n^{\eta-c_4}+n^{-c_8} ) ,\label{A.30}\hspace*{-35pt}
\end{equation}
where we used
\eqref{A.26} to derive the final identity. Therefore, \eqref{A.4}
holds with ${\lambda }_n=n^{1-c_9}$ for any $c_9\in(0,\max(c_4,c_8))$.
Hence we may use this value of ${\lambda }_n$ in \eqref{A.7}, establishing
that
%
\begin{equation}n\mo \sumjon\{g(X_j)-\hg(X_j)\}^2
=O_p (n^{-c} ) ,\label{A.31}
\end{equation}
with $c=\min(c_9,2a_1c_4)$,
where the estimator $\hbe$ used to define $\hg$ is obtained by
minimizing $S(\be)=\sumj(Y_j-\bY_j)^2$ [the first quantity in
\eqref{2.7}] over $\be\in\cB_1$. [We used \eqref{A.26} to simplify
the term in $h^{2a_1}$ in \eqref{A.7}.]

During the proof above we imposed on the class $\cB_1$ the
assumption that $\be\z\in\cB_1$ [see the discussion following
\eqref{A.5}], and also three conditions---\eqref{A.11},
\eqref{A.14} and \eqref{A.28}---on the size of the class. The latter three conditions hold if
%
\begin{equation}\#\cB_1=O \{\exp (n^{c_{10}} ) \}
,\label{A.32}
\end{equation}
provided $0<c_{10}<\min(1-c_2-c_5,1-a_4c_5,c_3)$. (Recall from
Step~6 that $m$ equals the integer part of $n^{c_3}$.) By choosing
$c_5$ smaller if necessary we can ensure that the upper bound here
is strictly positive, and so $c_{10}>0$.

Let $0<c_{11}<c_{10}$ and $c_{12}>0$, define $r=r(n)$ to be the
integer part of $n^{c_{11}}$, and let $D_3$ be as in \eqref{3.4}.
Let $r$ be as stipulated in \eqref{3.8}, and write $\cB_2$ for the
class of functions $\be=\sum_{1\leq k\leq r} b_k \psi_k$ such that
each $|b_k|\leq D_3$ for $1\leq k\leq r$. Let $\cB_3$ be the set of
elements of $\cB_2$ for which each $b_k$, for $1\leq k\leq r$, is an
integer multiple of $n^{-c_{12}}$. The number of elements of
$\cB_3$ is bounded above by a constant multiple of
%
\begin{equation}
 (2D_3 n^{c_{12}} )^r\leq\exp (\mathrm{const.}\ n^{c_{11}} \log
n ) =o \{\exp (n^{c_{10}} ) \} .\label{A.33}
\end{equation}
Put
$\cB_1=\cB_3\cup\{\be\z\}$. Then \eqref{A.32} follows from
\eqref{A.33}.

The following three properties hold: (a) The lattice on which
$\cB_3$ is based can be made arbitrarily fine in a polynomial sense,
by choosing $c_{12}$ sufficiently large; (b) $E\|X\|^\eta<\infty$
for some $\eta>0$ [see \eqref{3.3}]; and (c) $K$ has a bounded
derivative [see \eqref{3.7}]. Given $\be=\sum_{1\leq k\leq
r} b_k \psi_k\in\cB_2$, let $\be\approx=\sum_{1\leq k\leq
r} b_k\approx \psi_k$ be the element of $\cB_2$ defined by taking
$b_k\approx$ to be the lattice value nearest to $b_k$, for $1\leq
k\leq r$. Define $S(n)$ to equal the maximum, over $1\leq i,j\leq
n$, of $\|X_i-X_j\|$. Property (b) implies that
$S_n=O_p(n^{c_{13}})$ for some $c_{13}>0$. Using this property, and
(a) and (c), it can be proved, by taking $c_{12}$ sufficiently
large, that for any given $c_{14}>0$,
%
\begin{eqnarray}\label{A.34}
 && \sup_{\be\in\cB_2} \max_{1\leq i,j\leq n}
 |K\ij(\be)-K\ij (\be\approx ) |\nonumber\\
&& \qquad =O_p \Bigl\{S(n) h\mo \sup_{\be\in\cB_2}  \|\be-\be
\approx \|
 \Bigr\}\\
   && \qquad =O_p (n^{-c_{14}} ) .
\nonumber
\end{eqnarray}
From this result and the other properties of $K$ in \eqref{3.7} it
can be shown that \eqref{A.31} continues to hold if $\hbe$ in the
definition of $\hg$ is replaced by the minimizer of
$S(\be)=\sumj(Y_j-\bY_j)^2$ over $\be\in\cB_4=\cB_2\cup\{\be\z
\}$.

Call this result (R).

The desired result \eqref{3.9} follows from (R), except that the set
$\cB_4$ contains $\be\z$ as an unusual, adjoined element. Hence
there is, in theory, a possibility that $\hbe=\be\z$; this could not
happen if we were to restrict $\hbe$ to elements of $\cB_2$, as
required when defining the estimator $\hg$ in \eqref{3.9}. To
appreciate that this does not cause any difficulty, let
$\be^1=\sum_{1\leq k\leq r} b_k\z \psi_k$ denote the approximation
to $\be\z$ obtained by dropping all but the first $r$ terms in the
expansion $\be\z=\sum_{k\geq1} b_k\z \psi_k$. The argument
leading to \eqref{A.34} can be used to prove that, for $c_{14}>0$
chosen arbitrarily large, there exists a value of $B=B(c_{14})$, in
the second part of \eqref{3.4}, such that
\[\max_{1\leq i,j\leq
n}   |K\ij(\be\z)-K\ij(\be^1) |
=O_p \{S(n) h\mo  \|\be\z-\be^1 \| \}
=O_p (n^{-c_{14}} ) .
\]
Arguing as before, this leads to
the conclusion that $\be\z$ can be dropped from $\cB_4$ without
damaging result (R).
\end{appendix}

\section*{Acknowledgment}
We wish to thank two reviewers for very helpful comments and suggestions.


%

\printaddresses

\end{document}